\documentclass[11pt,oneside]{article}
\usepackage[square, numbers]{natbib}
\usepackage[a4paper, margin=1in]{geometry}
\usepackage[T1]{fontenc}
\usepackage{amsmath, amssymb, amsthm,caption}
\usepackage{xcolor}
\usepackage{mathtools}
\usepackage{amsfonts}
\usepackage{hyperref}
\usepackage{enumerate}
\usepackage{bm, bbm, cancel}
\usepackage{hhline}             
\usepackage[french, english]{babel}

\providecommand{\keywords}[1]
{
  \small	
  \textbf{\textit{Keywords---}} #1
}

\usepackage[]{todonotes}

\newcommand{\donghan}[2][]{\todo[color=red!20!white, size=\footnotesize, #1]{DK: #2}}

\newcommand{\R}{\mathbb{R}}
\newcommand{\N}{\mathbb{N}}

\newcommand{\T}{\mathbb{T}}

\numberwithin{equation}{section}

\newtheorem{theorem}{Theorem}[section] 
\newtheorem{lemma}[theorem]{Lemma}
\newtheorem{corollary}[theorem]{Corollary}

\newtheorem{proposition}[theorem]{Proposition}

\theoremstyle{definition}
\newtheorem{definition}[theorem]{Definition}

\newtheorem{remark}[theorem]{Remark}

\newtheorem{examplex}{Example}
\newenvironment{example}
{\pushQED{\qed}\examplex}
{\popQED\endexamplex}

\title{On isomorphism of the space of continuous functions with finite \texorpdfstring{$p$}-th variation along a partition sequence.}

\vspace{2mm}

\author{
	\textsc{Purba Das}
	\thanks{Department of Mathematics, King's College London, UK (E-mail: \it{purba.das@kcl.ac.uk})}
	\and
	\textsc{Donghan Kim} 
	\thanks{Department of Mathematical Sciences, KAIST, South Korea (E-mail: {\it kimdonghan@kaist.ac.kr})}
}

\begin{document}
\maketitle
\begin{abstract}
		We study the concept of (generalized) $p$-th variation of a real-valued continuous function along a general class of refining sequence of partitions. We show that the finiteness of the $p$-th variation of a given function is closely related to the finiteness of $\ell^p$-norm of the coefficients along a Schauder basis, similar to the fact that H\"older coefficient of the function is connected to $\ell^{\infty}$-norm of the Schauder coefficients. This result provides an isomorphism between the space of $\alpha$-H\"older continuous functions with finite (generalized) $p$-th variation along a given partition sequence and a subclass of infinite-dimensional matrices equipped with an appropriate norm, in the spirit of Ciesielski.

    \end{abstract}

    \smallskip
    
\keywords{\texorpdfstring{$p$}{TEXT}-th variation, H\"older regularity, Ciesielski's isomorphism, Schauder basis, Variation index, Refining partition sequences}
\tableofcontents

    
\section{Introduction}
In the seminal paper \cite{follmer1981}, F\"ollmer derived the pathwise It\^o's formula for a class of real functions with a finite quadratic variation. In particular, for a twice differentiable function $F$ and a one-dimensional continuous function $x$ with finite quadratic variation along a partition sequence $\pi = (\pi^n)_{n \in \mathbb{N}}$, the pathwise It\^o formula is given as
\begin{equation}   \label{eq : Follmer-Ito formula}
		F\big(x(t)\big) = F\big(x(0)\big) + \int_0^t F'\big(x(s)\big) d^{\pi} x(s) + \frac{1}{2} \int_0^t F''\big(x(s)\big) d[x]_{\pi}(s).
\end{equation}
    Here, the first integral is defined as a left Riemann sum
\begin{equation*}
	\int_0^t F'\big(x(s)\big) d^{\pi}x(s) := \lim_{n \to \infty} \sum_{\pi^n \ni t^n_j \le t} F'\big(x(t^n_j)\big) \big( x(t^n_{j+1}) - x(t^n_{j}) \big),
\end{equation*}
    
and the integrator $[x]_{\pi}(\cdot)$ of the second integral is the quadratic variation of $x$ along the partition sequence $\pi$, defined as the following uniform limit in $t$:
    \begin{equation}    \label{def: quadratic variation}
		[x]_{\pi^n}(t) := \sum_{\pi^n \ni t^n_j \le t} \big \vert x(t^n_{j + 1}) - x(t^n_j)\big\vert^2 \xlongrightarrow{n \rightarrow \infty} [x]_\pi(t).
    \end{equation}

    This pathwise It\^o's formula has been generalized in several aspects \cite{ananova2017, chiu2018, CF2010, perkowski2019, davis2018, Kim:localtime, Rafal2021}. Among these, Cont and Perkowski \cite{perkowski2019} defined the notion of $p$-th variation of continuous functions along $\pi$ by raising the exponent in \eqref{def: quadratic variation} to any even integers $p \in 2\mathbb{N}$, and derived high-order pathwise change-of-variable formula; more recently, Cont and Jin \cite{ruhong2022} developed fractional pathwise It\^o formula for functions with $p$-th variation for any $p > 1$, with a fractional It\^o remainder term. These pathwise calculus formulae, including F\"ollmer's original one \eqref{eq : Follmer-Ito formula}, require the continuous function $x$ to have finite $p$-th variation along $\pi$. In other words, the existence of the limit
\begin{equation}    \label{def: p-th variation}
[x]^{(p)}_{\pi^n}(t) := \sum_{\pi^n \ni t^n_j \le t} \big \vert x(t^n_{j + 1}) - x(t^n_j)\big\vert^p \xlongrightarrow{n \rightarrow \infty} [x]^{(p)}_\pi(t)
\end{equation}
is the crucial assumption when applying these formulae. It is then natural to study a class $V^{p}_{\pi}$ of functions $x$ such that the limit \eqref{def: p-th variation} exists for a fixed partition sequence $\pi$ and $p > 1$. 
	
In this regard, Schied \cite{schied2016} showed that the space $V^{p}_{\pi}$ is not a vector space by constructing an example of two continuous functions $x$ and $y$ on $[0, 1]$ such that $[x]^{(2)}_{\mathbb{T}}$ and $[y]^{(2)}_{\mathbb{T}}$ exist, but $[x+y]^{(2)}_{\mathbb{T}}$ does not exist, along the dyadic partition sequence $\mathbb{T} = (\mathbb{T}^n)_{n \in \mathbb{N}}$ with $\mathbb{T}^n := \{k2^{-n} : k = 0, 1, \cdots, 2^n\}$. These two functions $x$ and $y$ belong to a class of so-called generalized Takagi functions, constructed via the Schauder representation of continuous functions. From the Schauder representation of $x$ and $y$ along $\mathbb{T}$, one can obtain explicit expressions of both terms in the following strict inequality to show that $[x+y]^{(2)}_{\mathbb{T}}$ does not exist:
	\begin{equation*}
		\liminf_{n \to \infty} \, [x+y]^{(2)}_{\mathbb{T}^n}(t) < \limsup_{n \to \infty} \, [x+y]^{(2)}_{\mathbb{T}^n}(t).
	\end{equation*}
    Since Schied's example implies that requiring the existence of the limit \eqref{def: p-th variation} restricts the function space $V^p_{\pi}$ too much, in this paper we study a larger space $\mathcal{X}^p_{\pi} \supset V^p_\pi$ of functions $x$ that satisfy 
    \begin{equation}	\label{con : limsup finite}
		\limsup_{n \to \infty} \, [x]^{(p)}_{\pi^n}(t) = \limsup_{n \to \infty} \sum_{\pi^n \ni t^n_j \le t} \big \vert x(t^n_{j + 1}) - x(t^n_j)\big\vert^p < \infty,
    \end{equation}
    but does not require the limit to exist. With an appropriate norm, we prove that the space $\mathcal{X}^p_{\pi}$ is a Banach space (see definition \eqref{def : x p pi space} and Proposition \ref{prop : Banach space} below).
	
    Even though we may not apply the aforementioned pathwise change-of-variable formulae to every function in $\mathcal{X}^p_{\pi}$, we shall study the Banach space $\mathcal{X}^p_{\pi}$, instead of $V^p_{\pi}$, because the notion of variation index, i.e., the infimum number $p \ge 1$ such that the condition \eqref{con : limsup finite} holds (see Definition \ref{def. variation index} below), can be used for measuring `roughness' of a given function (or a path of a stochastic process) \cite{fake_fBM,das2022theory}. It is well known that (almost every path of) a fractional Brownian motion (fBM) $B^H$ with Hurst index $H \in (0, 1)$, has H\"older exponent equal to $H-$, whereas its variation index along `reasonable' partition sequences (e.g., dyadic partition sequence $\mathbb{T}$) is equal to $1/H$. These facts are closely related to the self-similarity property of fBMs, but it is generally not true for general continuous functions that the reciprocal of the variation index is equal to (the supremum of) H\"older exponent. In a recent work \cite{fake_fBM}, a specific example of $(1/4)$-H\"older continuous function with variation index along the dyadic partition sequence equal to $2$ is constructed, thus, the variation index should be considered as an alternative way of measuring function's roughness.

With the help of Schauder representation along a general class of partition sequences, our main result provides a necessary and sufficient condition for elements of the Banach space $\mathcal{X}^p_{\pi}$, in terms of their Schauder coefficients (see Theorem \ref{thm: p-th var general2}). More specifically, the condition \eqref{con : limsup finite} is equivalent to the $\ell^{\infty}$-finiteness of the sequence composed of $\ell^p$-norm of Schauder coefficients of functions along each partition $\pi^n$, scaled by a $(p/2)$-power of the mesh size of $\pi^n$.
	
    When the Schauder coefficients of functions are arranged in an infinite dimensional matrix, this result gives rise to an isomorphism between the space of $\alpha$-H\"older continuous functions with finite (generalized) $p$-th variation along a partition sequence $\pi$ and a subspace of infinite-dimensional matrices with an appropriate matrix norm (see Theorem \ref{thm : isomorphism}). Our isomorphism result reminds that of Ciesielski's in 1960 \cite{Ciesielski:isomorphism}, between the space of $\alpha$-H\"older continuous functions and the space of bounded real sequences, using Schauder representation along the dyadic partition sequence $\mathbb{T}$, which has been generalized recently by \cite{fake_fBM} along a wider class of partition sequences.

    {\it Preview:} This paper is organized as follows. Section \ref{sec: variation index and the Banach space} introduces the notion of variation index and defines the Banach space $\mathcal{X}^p_{\pi}$. Section \ref{sec: Schauder} provides some notations and reviews preliminary results regarding Schauder representation of continuous functions. Section \ref{sec: main result} states and proves our main result, the characterization of generalized $p$-th variation in terms of Schauder coefficients along a given partition sequence. Section \ref{sec: isomorphism} includes the isomorphism, as an important consequence of the result. 

\bigskip
	
\section{Variation index and the Banach space \texorpdfstring{$\mathcal{X}^p_{\pi}$}{TEXT}} \label{sec: variation index and the Banach space}
	
\medskip
	
\subsection{\texorpdfstring{$p$}{TEXT}-th variation and variation index}    \label{subsec: variation index}
	
First, we introduce some relevant notations and definitions for partition sequences. For a fixed $T>0$, we shall consider a (deterministic) sequence of partitions $\pi=(\pi^n)_{n \ge 0}$ of $[0,T]$
\begin{equation*}
\pi^n = \left(0=t^n_0<t^n_1<t^n_2<\cdots<t^n_{N(\pi^n)}=T\right),
\end{equation*}
where we denote $N(\pi^n)$ the number of intervals in the partition $\pi^n$. By convention, $\pi^0 = \{0, T\}$. For example, the dyadic partition sequence, denoted by $\mathbb{T} \equiv \pi$, contains partition points $t^n_k = kT/2^n$ for $n \in \mathbb{N}\cup \{0\}$, $k = 0, \cdots, 2^n$.
	
\begin{definition} [Refining sequence of partitions]	\label{Def : refining partition}
A sequence of partitions $\pi=(\pi^n)_{n \ge 0}$ is said to be \textit{refining (or nested)}, if $t\in \pi^m$ implies $t\in \cap_{n\geq m} \pi^n$ for every $m \in \mathbb{N}$. In particular, we have $\pi^0\subseteq\pi^1\subseteq \pi^2 \subseteq \cdots$. 
\end{definition}
	
For a partition sequence $\pi=(\pi^n)_{n \ge 0}$, we write
	\begin{equation}	\label{def : mesh size}
		\underline{\pi^n} := \inf_{i = 0, \cdots, N(\pi^n)-1} \vert t^n_{i+1} - t^n_i \vert, \qquad \qquad
		|\pi^n| := \sup_{i = 0, \cdots, N(\pi^n)-1} \vert t^n_{i+1} - t^n_i \vert,
	\end{equation}
	the size of the smallest and the largest interval of $\pi^n$, respectively. In the following, we denote $\Pi([0, T])$ the collection of all refining partition sequences $\pi$ of $[0, T]$ with vanishing mesh, i.e., $|\pi^n| \rightarrow 0$ as $n \rightarrow \infty$. 
	
	Let us denote $C^0([0, T])$ the space of real-valued continuous functions defined on $[0, T]$. In this subsection, we fix a partition sequence $\pi = (\pi^n)_{n \ge 0} \in \Pi([0, T])$ and $x \in C^0([0, T])$. For $p \ge 1$, we denote 
	\begin{equation}    \label{def: discrete p-th variation}
		[x]_{\pi^n}^{(p)}(t) := \sum_{\pi^n \ni t^n_j \le t} \big \vert x(t^n_{j + 1}) - x(t^n_j)\big\vert^p 
	\end{equation}
	the $p$-th variation of $x$ along a partition $\pi^n$ for each level $n \in \mathbb{N}$.

\begin{remark}	\label{rem: lim of p-th variation}
If there exists a continuous, non-decreasing function $[x]^{(p)}_{\pi}$ such that 
\begin{equation}	\label{eq : uniform limit of p-th variation}
	\lim_{n \to \infty} [x]_{\pi^n}^{(p)}(t) = [x]^{(p)}_{\pi}(t), \qquad \forall \, t \in [0, T],
\end{equation}
then we say $x$ admits finite $p$-th variation along $\pi$, and the above convergence is uniform in $t$ (\cite[Definition 1.1 and Lemma 1.3]{perkowski2019}). We write $V^{p}_{\pi}$ the space of such functions $x$ admitting finite $p$-th variation along $\pi$. In the particular case of $p = 2$ (then $V^2_{\pi}$ is often denoted as $Q_{\pi}$) and $\pi$ given as the dyadic partition sequence $\T$, it is shown in \cite[Proposition 2.7]{schied2016} that $V^{2}_{\T}$ is not a vector space. Note that this notion of the $p$-th variation \eqref{eq : uniform limit of p-th variation} is different from the $p$-variation commonly used in Lyons’ rough path theory, as the latter is defined as the supremum of $p$-th variation over all possible partitions of $[0, T]$.
	\end{remark}
	
	Even though the $p$-th variation of $x$ along a given sequence $\pi$ defined in Remark~\ref{rem: lim of p-th variation} may not exist, one can always define its variation index along $\pi$ as follows. Note that a recent paper \cite{10.1214/24-AAP2135} introduces a similar
concept to the variation index, referred to instead as the ‘roughness exponent’.
	
	\begin{definition} [Variation index along a partition sequence, Definition 2.3 of \cite{das2022theory}]  \label{def. variation index}
		The \textit{variation index} of $x \in C^0([0, T])$ along $\pi \in \Pi([0, T])$ is defined as
		\begin{equation}    \label{def. alter variation index}
			p^{\pi}(x) := \inf \big\{ p\geq 1 : \limsup_{n \to \infty} \, [x]_{\pi^n}^{(p)}(T) < \infty \big\}.
		\end{equation}
	\end{definition}

\noindent	Thanks to the continuity of $x$, it is straightforward to show
\begin{equation}    \label{eq : limsup q variation}
\limsup_{n \to \infty} \,[x]_{\pi^n}^{(q)}(T) =
	\begin{cases}
		\, 0, \qquad & q > p^{\pi}(x),
		\\
		\infty, \qquad & q < p^{\pi}(x),
	\end{cases}
\end{equation}
Therefore, the definition \eqref{def. alter variation index} can be re-formulated as 
\begin{equation*}
p^\pi(x) = \inf \big\{p\geq 1 \, : \, \limsup_{n \to \infty} \,[x]_{\pi^n}^{(p)}(T) = 0 \big\}.
\end{equation*}
Moreover, since $\limsup_{n \to \infty} \, [x]_{\pi^n}^{(p)}(T) < \infty$ if and only if $\sup_{n \in \N} \, [x]_{\pi^n}^{(p)}(T) < \infty$, we also have
\begin{equation}    \label{def. alter variation index2}
	p^{\pi}(x) = \inf \big\{ p\geq 1 : \sup_{n \in \N} \, [x]_{\pi^n}^{(p)}(T) < \infty \big\}.
\end{equation}	
Now that the quantity $[x]_{\pi^n}^{(p)}(t)$ in \eqref{def: discrete p-th variation} can be recognized as the $p$-th power of $\ell^p$-norm of the real sequence $\{x(t^n_{j+1})-x(t^n_j)\}_{t^n_j \in \pi^n, \, t^n_j \le t}$, we provide the following definition.
	
\begin{definition}\label{Def : p norm and Xp space}
For $x \in C^0([0, T])$, $p \ge 1$, and $\pi \in \Pi([0, T])$, we denote
\begin{equation*}
\Vert x \Vert^{(p)}_{\pi} := |x(0)| + \sup_{n \in \N} \, \Big([x]_{\pi^n}^{(p)} (T)\Big)^{\frac{1}{p}}
\end{equation*}
and consider the subspace of $C^0([0, T])$:
\begin{equation}    \label{def : x p pi space}
\mathcal X^{p}_{\pi} := \{x \in C^0([0, T])\, :  \Vert x \Vert^{(p)}_{\pi} < \infty\}.
\end{equation}
We say $\mathcal X^{p}_{\pi}$ is the class of continuous functions with finite \textit{(generalized) $p$-th variation} along $\pi$.
\end{definition}

The space $\mathcal X^{p}_{\pi}$ turns out to be a Banach space, in contrast to the space $V^p_\pi$. 
\begin{proposition}	\label{prop : Banach space}
The mapping $\mathcal X^{p}_{\pi} \ni x \mapsto \Vert x \Vert^{(p)}_{\pi}$ is a norm, and the space $(\mathcal X^{p}_{\pi}, \, \Vert \cdot \Vert^{(p)}_{\pi})$ is a Banach space.
\end{proposition}

\begin{proof}
    We first prove that the mapping is a norm. For any scalar $r$, the identity $\Vert r x \Vert^{(p)}_{\pi} = |r| \Vert x\Vert^{(p)}_{\pi}$ is straightforward. Thanks to Minkowski's inequality, it is also easy to prove the subadditive property (triangle inequality). These imply, in particular, that $\mathcal X^{p}_{\pi}$ is a vector space. Finally, if $\Vert x \Vert^{(p)}_{\pi} = 0$, then $x$ has zero value on every partition point $t^n_j$ of $\pi$ for all $j, n$. Since $|\pi^n| \rightarrow 0$ as $n \to \infty$, the set $P := \bigcup_{n \in \N} \pi^n$ of all partition points of $\pi$ is dense in $[0,T]$, and the continuity of $x$ with $x(0) = 0$ concludes $x \equiv 0$. This shows that $\Vert x \Vert^{(p)}_{\pi}$ is a norm.
		
    To prove the space $\mathcal{X}^p_{\pi}$ is a Banach space, we fix a Cauchy sequence $(x_{\ell})_{\ell \in \N}$ of $\mathcal X^{p}_{\pi}$, i.e., for any $\epsilon > 0$, there exists $N \in \mathbb{N}$ such that $\Vert x_k - x_m \Vert^{(p)}_{\pi} < \epsilon$ for all $k, m \ge N$. In particular, for every $k, m \ge N$, we have $|x_k(0)-x_m(0)| < \epsilon$ and 
		\begin{equation}    \label{eq : x, m bound}
			[x_k - x_m]^{(p)}_{\pi^n}(T) = \sum_{t^n_j \in \pi^n} \Big \vert \big( x_k(t^n_{j + 1}) - x_m(t^n_{j + 1}) \big) - \big( x_k(t^n_j) - x_m(t^n_j) \big) \Big\vert^p < \epsilon^p
		\end{equation}
		holds for each $n \in \N$. Since $\{x_{\ell}(0)\}_{\ell \in \N}$ is a real Cauchy sequence, its limit $\lim_{\ell \to \infty} x_{\ell}(0) = \tilde{x}(0)$ exists. Moreover, we fix an arbitrary $n \in \N$, then for all indices $j$ such that $t^n_j$ belongs to $\pi^n$, we have
		\begin{align*}
			\Big \vert \big( x_k(t^n_{j + 1}) - x_k(t^n_{j}) \big) - \big( x_m(t^n_{j+1}) - x_m(t^n_j) \big) \Big\vert^p \\= \Big \vert \big( x_k(t^n_{j + 1}) - x_m(t^n_{j + 1}) \big) - \big( x_k(t^n_j) - x_m(t^n_j) \big) \Big\vert^p < \epsilon^p
		\end{align*}
		for every $k, m \ge N$, in other words, $\big(x_k(t^n_{j + 1}) - x_k(t^n_{j}) \big)_{k \in \mathbb{N}}$ is a Cauchy sequence in $\mathbb{R}$ for each $j$. Again by the completeness of $\mathbb{R}$, the limit $d(t^n_j) := \lim_{k \to \infty} \big( x_k(t^n_{j+1}) - x_k(t^n_{j}) \big) \in \R$ exists for each index $j$ and $n \in \N$.
		
		Let us recall the set $P = \bigcup_{n \in \N} \pi^n$ of all partition points of $\pi$, and define a function $\tilde{x}$ on $P$
		\begin{equation*}
			\tilde{x}(t^n_j) = \tilde{x}(0) + \sum_{i=0}^{j-1} d(t^n_i), \qquad \text{for every } t^n_j \in \pi^n \text{ and } n \in \N.
		\end{equation*}
		Since $P$ is a dense subset of $[0, T]$ and a function defined on a dense set can be extended to a continuous function, there exists $x \in C^0([0, T])$ such that $x(t^n_j) = \tilde{x}(t^n_j)$ holds for all points $t^n_j$ of $P$. Furthermore, we have $x(0) = \tilde{x}(0) = \lim_{k \to \infty} x_k(0)$ as well as
		\begin{equation*}
			x(t^n_{j+1}) - x(t^n_{j}) = \tilde{x}(t^n_{j+1}) - \tilde{x}(t^n_{j})
			= d(t^n_j) = \lim_{k \to \infty} \big( x_k(t^n_{j+1}) - x_k(t^n_{j}) \big),
		\end{equation*}
		thus $x(t^n_j) = \lim_{k \to \infty} x_k(t^n_j)$ for each $t^n_j \in P$.
		
		Sending $m \to \infty$ in \eqref{eq : x, m bound}, we have for each $n \in \N$
		\begin{equation}    \label{ineq : x_k and x}
			\sum_{t^n_j \in \pi^n} \Big \vert \big( x_k(t^n_{j + 1}) - x(t^n_{j + 1}) \big) - \big( x_k(t^n_j) - x(t^n_j) \big) \Big\vert^p < \epsilon^p, \qquad \text{for } k \ge N.
		\end{equation}
		Minkowski's inequality now yields for each $n \in \N$
		\begin{align*}
			\bigg( \sum_{t^n_j \in \pi^n} \Big \vert x(t^n_{j + 1}) - x(t^n_{j}) \Big\vert^p \bigg)^{\frac{1}{p}}
			&\le \bigg( \sum_{t^n_j \in \pi^n} \Big \vert \big( x_k(t^n_{j + 1}) - x(t^n_{j + 1}) \big) - \big( x_k(t^n_j) - x(t^n_j) \big) \Big\vert^p \bigg)^{\frac{1}{p}}
			\\
			& \qquad + \bigg( \sum_{t^n_j \in \pi^n} \Big \vert x_k(t^n_{j + 1}) - x_k(t^n_{j}) \Big\vert^p \bigg)^{\frac{1}{p}}
			\\
			&\le \epsilon + \Vert x_k \Vert^{(p)}_{\pi} < \infty, \qquad \text{for } k \ge N,
		\end{align*}
		and this proves $x \in \mathcal{X}^{p}_{\pi}$. Furthermore, the inequality \eqref{ineq : x_k and x} implies $\Vert x_k - x \Vert^{(p)}_{\pi} < \epsilon$ for all large enough numbers $k$. This concludes that the Cauchy sequence $(x_{\ell})_{\ell \in \N}$ converges to $x$ in $\Vert \cdot \Vert^{(p)}_{\pi}$ norm.
	\end{proof}
	
	In line with Proposition \ref{prop : Banach space}, it is well-known that the space $(C^{0, \alpha}([0, T]), \, \Vert \cdot \Vert_{C^{0, \alpha}})$ of $\alpha$-H\"older continuous functions, is also a Banach space. We next note the inclusion
	\begin{equation}	\label{inclusion}
		\mathcal{X}^{p}_{\pi} \subset \mathcal{X}^{q}_{\pi}, \qquad \text{for } 1 \le p \le q < \infty,
	\end{equation}
	due to the straightforward inequality $([x]^{(q)}_{\pi^n}(T))^{\frac{1}{q}} \le ([x]^{(p)}_{\pi^n}(T))^{\frac{1}{p}}$ for every $n \ge 0$. We conclude this subsection with the following property that adding a function with vanishing $p$-th variation does not affect the variation index.
\begin{lemma}   \label{lem: adding zero p-th variation}
    For $x, y \in C^0([0, T])$, $p \ge 1$, $t \in [0, T]$, and $\pi \in \Pi([0, T])$, suppose that 
    \begin{equation*}
        \limsup_{n \to \infty} \, [y]^{(p)}_{\pi^n}(t) = 0  
    \end{equation*}
    holds. Then, we have
    \begin{equation*}
        \limsup_{n \to \infty} \, [x]^{(p)}_{\pi^n}(t) < \infty \qquad \text{if and only if} \qquad \limsup_{n \to \infty} \, [x+y]^{(p)}_{\pi^n}(t) < \infty,
    \end{equation*}
    therefore $p^{\pi}(x) = p^{\pi}(x+y)$. In particular, the identity $[x]^{(p)}_{\pi}(t) = [x+y]^{(p)}_{\pi}(t)$ holds, provided that the limit $[x]^{(p)}_{\pi}(t)$ exists in the sense of Remark \ref{rem: lim of p-th variation}.
\end{lemma}
\begin{proof}
    Applying Minkowski's inequality twice yields
    \begin{equation*}
        \big( [x]^{(p)}_{\pi^n}(t) \big)^{\frac{1}{p}} - \big( [y]^{(p)}_{\pi^n}(t) \big)^{\frac{1}{p}} \le \big( [x+y]^{(p)}_{\pi^n}(t) \big)^{\frac{1}{p}} \le \big( [x]^{(p)}_{\pi^n}(t) \big)^{\frac{1}{p}} + \big( [y]^{(p)}_{\pi^n}(t) \big)^{\frac{1}{p}}.
    \end{equation*}
    Taking $\limsup$ or $\lim$ respectively gives the result.
\end{proof}

	\medskip
	
	\subsection{Variation index along different partition sequences}
	
	A continuous function $x$ can have different $p$-th variations, $[x]^{(p)}_{\pi}$ and $[x]^{(p)}_{\rho}$, along two different refining partition sequences $\pi$ and $\rho$. In this subsection, we study the variation index of $x$ along different partition sequences. We first introduce Proposition \ref{prop : partition with q variation zero}, inspired by Freedman \cite{freedman}, whose proof needs a preliminary result.
	
	\begin{lemma}   \label{lem : partition for any q}
		For any given numbers $q>1$, $\epsilon > 0$, and $x \in C^0([0, T])$, there exists a finite set $\pi = \{0 = t_0< t_1< \cdots< t_m = T\}$ in $[0, T]$ such that the $q$-th variation of $x$ along $\pi$ is less than $\epsilon$, i.e.,
		\begin{equation*}
			[x]^{(q)}_{\pi}(T) = \sum_{j = 0}^{m-1} \Big \vert x(t_{j + 1}) - x(t_j)\Big\vert^q < \epsilon.
		\end{equation*} 
	\end{lemma}
	
	\begin{proof}
		If $x(0) = x(T)$, then we just take $\pi = \{0, T\}$. Thus, we suppose that $x(T) > x(0)$; the other case $x(T) < x(0)$ can be handled by applying the same argument to $y(t) = x(T-t)$.
		
		We assume without loss of generality that $x(0) = 0$, $T=1$, and $x(T) = 1$. For given $q > 1$ and $\epsilon > 0$, we choose $N \in \mathbb{N}$ large enough so that $N^{1-q} < \epsilon$, and define $t^N_{j} := \min\{t \ge 0 : x(t) = j/N\}$ for $j = 0, \cdots, N$. Let $\pi = \{t^N_{0}, \cdots, t^N_{N}\}$ if $t^N_{N} = 1$, or $\pi = \{t^N_{0}, \cdots, t^N_{N}, 1\}$ otherwise. Now it is simple to check $[x]^{(q)}_{\pi}(1) = N^{1-q} < \epsilon$.
	\end{proof}
	
\begin{proposition} \label{prop : partition with q variation zero}
For any $x \in C^0([0, T])$, we have
		\begin{equation*}
			\inf \big\{ p^{\pi}(x) : \pi \in \Pi([0, T]) \big\} = 1.
		\end{equation*}
	\end{proposition}
	
	\begin{proof}
		Let us fix $x \in C^0([0, T])$. For any $q > 1$, we shall show that there exists a sequence $\pi = (\pi^n)_{n \ge 0} \in \Pi([0, T])$ satisfying
		\begin{equation}    \label{eq : decreasing partition to zero}
			[x]^{(q)}_{\pi}(T) = \limsup_{n \to \infty} \, [x]^{(q)}_{\pi^n}(T) = 0.
		\end{equation}
		Then, the identity \eqref{eq : decreasing partition to zero}, together with \eqref{eq : limsup q variation}, implies that for any $q > 1$ there exists $\pi \in \Pi([0, T])$ satisfying $p^{\pi}(x) \le q$, which in turn proves the result.
		
		We choose a decreasing real sequence $\epsilon_n \downarrow 0$, and set $\pi^0 = \{0, T\}$. We shall inductively define $\pi^n$ for each $n \ge 0$. Suppose $\pi^n$ is defined, and let $\rho^{n+1}$ be a partition of $[0, T]$ satisfying $\pi^n \subset \rho^{n+1}$ and $|\rho^{n+1}| \le \epsilon_{n+1}$. Suppose that $\rho^{n+1}$ has $m+1$ points, dividing $[0, T]$ into $m$ subintervals. From Lemma \ref{lem : partition for any q}, we construct a partition $\pi^{n+1}$ of $[0, T]$ with $\rho^{n+1} \subset \pi^{n+1}$, such that for each pair $t^{\rho^{n+1}}_j, t^{\rho^{n+1}}_{j+1}$ of consecutive points of $\rho^{n+1}$ we have 
		\begin{equation*}
			[x]^{(q)}_{\nu^{n+1}_j} \le \frac{\epsilon_{n+1}}{m},
		\end{equation*}
		where $\nu^{n+1}_j := \pi^{n+1} \cap [t^{\rho^{n+1}}_{j}, t^{\rho^{n+1}}_{j+1}]$ and $[x]^{(q)}_{\nu^{n+1}_j}$ is the $q$-th variation along $\nu^{n+1}_j$ on the interval $[t^{\rho^{n+1}}_{j}, t^{\rho^{n+1}}_{j+1}]$. Then, we obtain $[x]^{(q)}_{\pi^{n+1}}(T) \le \epsilon_{n+1}$ and $|\pi^{n+1}| \le |\rho^{n+1}| \le \epsilon_{n+1}$, therefore, $\pi = (\pi^n)$ satisfies condition \eqref{eq : decreasing partition to zero}.
	\end{proof}
The partition sequence $\pi$, constructed in the proof of Proposition \ref{prop : partition with q variation zero}, is a Lebesgue-type partition of $x$.
A recent paper \cite{rafal2025} presents an example of a continuous function, known as the Peano curve, which exhibits
different values of quadratic variation along different sequences of Lebesgue-type partitions.
    
	On the other hand, it is well known \cite{FrizHairer}  that an $\alpha$-H\"older continuous function $x \in C^{0, \alpha}([0, T])$ has finite $(\frac{1}{\alpha})$-variation, i.e., $\Vert x \Vert_{\frac{1}{\alpha}-var} < \infty$, with
	\begin{equation*}
		\Vert x \Vert_{p-var} := \bigg( \sup_{\rho} \sum_{t_j, t_{j+1} \in \rho} \big|x(t_{j+1}) - x(t_j) \big|^p \bigg)^{\frac{1}{p}},
	\end{equation*}
	where the supremum is taken over all partitions $\rho$ of $[0, T]$. This implies that for a given refining partition sequence $\pi \in \Pi([0, T])$ with vanishing mesh, the variation index $p^{\pi}(x)$ of $x \in C^{0, \alpha}([0, T])$ should be bounded above by the reciprocal of its H\"older exponent $\alpha$ (see Lemma~4.3 of \cite{fake_fBM} for the proof), namely
    \begin{equation*}
		p^{\pi}(x) \le \frac{1}{\alpha}.
    \end{equation*}

    We formalize the above arguments into the following theorem.
    \begin{theorem} \label{thm : variation index bounds}
		For any $x \in C^0([0, T])$, we have
		\begin{equation*}
			\inf \big\{p^{\pi}(x) : \pi \in \Pi([0, T]) \big\} = 1.
		\end{equation*}
		Moreover, for any $x \in C^{0, \alpha}([0, T])$, we have
		\begin{equation}  \label{ineq : sup variation index}
			\sup \big\{p^{\pi}(x) : \pi \in \Pi([0, T]) \big\} \le \frac{1}{\alpha}.
		\end{equation}
\end{theorem}

This result implies that an $\alpha$-H\"older continuous function $x$ can have any variation index $p^{\pi}(x)$ between $1$ and $1/\alpha$, along a given partition sequence $\pi \in \Pi([0, T])$. Moreover, the inclusion \eqref{inclusion} shows that $x \in \mathcal{X}^{q}_{\pi}$ for any $q > p^{\pi}(x)$.

\begin{example}
    The inequality \eqref{ineq : sup variation index} can be strict. Consider the increasing function $y(t) = \sqrt{t}$ defined on $[0, 1]$, which is $\frac{1}{2}$-H\"older continuous. The function $y$ has finite $1$-variation along any partition sequence $\pi$, thus $p^{\pi}(y)=1$, as it is an increasing function.
\end{example}

\begin{example}
A uniformly continuous function $z$ defined on $[0, \frac{1}{2}]$
\begin{equation*}
z(t) = \begin{cases}
\frac{1}{\log t}, \qquad & t \in (0, \frac{1}{2}],\\
~0, \qquad & t = 0,
\end{cases}
\end{equation*}
is not $\alpha$-H\"older continuous for any $\alpha > 0$. However, it is a decreasing function on the compact support, thus of bounded variation. As in the previous example, $p^{\pi}(z) = 1$ for every $\pi \in \Pi([0, \frac{1}{2}])$, which implies the left-hand side of \eqref{ineq : sup variation index} for $z$ is $1$.
\end{example}

In what follows, we shall characterize conditions for $x$ to belong to the Banach space $\mathcal{X}^{p}_{\pi}$, in terms of the Schauder coefficients of $x$ along $\pi$.

\bigskip
	
\section{Schauder representation along a general class of partition sequences}  \label{sec: Schauder}
	
In this section, we provide several definitions and preliminary results, mostly taken from \cite{das2021, das2020}, regarding Schauder representation of continuous functions along a general class of partition sequences. This type of representation was originally introduced by Schauder \cite{schauder1927}. After that, we shall provide our results in the next sections.
	
	\medskip
	
	\subsection{Properties of partition sequence}
	
	Let us recall Definition \ref{Def : refining partition} and the notations \eqref{def : mesh size}. We introduce a subclass of refining sequence of partitions with a `finite branching' property at every level $n \in \mathbb{N}$.
	
	\begin{definition} [Finitely refining sequence of partitions]    \label{def.finite.refining}
		A sequence of partitions $\pi=(\pi^n)_{n \ge 0}$ in $\Pi([0, T])$ is said to be \textit{finitely refining}, if there exists a positive integer $M$ such that the number of partition points of $\pi^{n+1}$ within any two consecutive partition points of $\pi^n$ is always bounded above by $M$, irrespective of $n \ge 0$. In particular, we have $\sup_{n \ge 0} \frac{N(\pi^n)}{M^n}\leq 1$.
	\end{definition}
	
	The following definition provides a condition that the ratio of the biggest step size to the smallest step size at each level is bounded.
	
	\begin{definition} [Balanced sequence of partitions] \label{def.balance}
		A sequence of partitions $\pi=(\pi^n)_{n \ge 0}$ is said to be \textit{balanced}, if there exists a constant $c>1$ such that 
		\begin{equation}    \label{eq.balance}
			\frac {|\pi^n|}{\underline{\pi^n}}\leq c
		\end{equation}
		holds for every $n \in \mathbb{N}$.
	\end{definition}
	
	We now give two conditions of refining partition sequences involving the biggest step sizes of two consecutive levels.
	
	\begin{definition} [Complete refining sequence of partitions]   \label{def.complete refining}
		A finitely refining sequence of partitions $\pi=(\pi^n)_{n \ge 0}$ is said to be \textit{complete refining}, if there exist positive constants $a$ and $b$ such that
		\begin{equation}    \label{def: complete refining}
			1+a \leq \frac{|\pi^n|}{|\pi^{n+1}|} \leq b
		\end{equation}
		holds for every $n \in \mathbb{N}$.
	\end{definition}
	
	\begin{definition} [Convergent refining sequence of partitions]		\label{def.conv.ref}
		A complete refining sequence of partitions is said to be \textit{convergent refining}, if the following limit exists: 
		\begin{equation}	\label{def: convergent refining}
			\lim_{n\to \infty} \frac{|\pi^n|}{|\pi^{n+1}|} = r \in (1, \infty). 
		\end{equation}
	\end{definition}
		
\begin{remark} [Notation]	\label{rem: constants}
Throughout this paper, we shall use the same symbols $M, c, a, b$, and $r$ to refer to the constants that appeared in Definitions \ref{def.finite.refining} - \ref{def.conv.ref}.
	\end{remark}
	
	\medskip
	
	\subsection{Generalized Haar and Schauder basis}	\label{subsec: Haar and Schauder}
	
	This subsection recalls some relevant definitions of generalized Haar and Schauder functions, which were introduced in \cite{das2021}.
	
	Let us fix $\pi \in \Pi([0, T])$ and denote $p(n,k) := \inf \{j \geq 0 : t^{n+1}_j \geq t^n_k \}$. Since $\pi$ is refining, we have the following inequality for every $k = 0, \cdots, N(\pi^n)-1$
	\begin{equation}    \label{EqFor_p}
		0\leq t^{n}_{k}= t^{n+1}_{p(n,k)}<t^{n+1}_{p(n,k)+1}<\cdots <t^{n+1}_{p(n,k+1)}= t^{n}_{k+1}\leq T.    
	\end{equation}
	With the notation $\Delta^{n}_{i, j} := t^n_j - t^n_i$, we now define the generalized Haar basis associated with $\pi$.
	
	\begin{definition} [Generalized Haar basis]	\label{Def : Haar}
		The \textit{generalized Haar basis} associated with a finitely refining sequence $\pi=(\pi^n)_{n \ge 0}$ of partitions is a collection of piecewise constant functions $\{\psi^{\pi}_{m,k,i} \, : \, m=0,1,\cdots, ~ k=0,\cdots,N(\pi^m)-1, ~ i = 1,\cdots, p(m,k+1)-p(m,k)\}$ defined as follows:
\begin{equation}	\label{haar_basis}
\psi^{\pi}_{m,k,i}(t)= 
	\begin{cases}
	\qquad \qquad \qquad \qquad 0, &\quad\text{if } t\notin \left[t^{m+1}_{p(m,k)},t_{p(m,k)+i}^{m+1}\right)
\\
	\quad \left( \frac{\Delta^{m+1}_{p(m, k)+i-1, p(m, k)+i}}{\Delta^{m+1}_{p(m, k), p(m, k)+i-1}} \times \frac{1}{\Delta^{m+1}_{p(m, k), p(m, k)+i}} \right)^{\frac{1}{2}}, &\quad\text{if } t\in\left[t_{p(m,k)}^{m+1},t_{p(m,k)+i-1}^{m+1}\right)
\\
	-\left( \frac{\Delta^{m+1}_{p(m, k), p(m, k)+i-1}}{\Delta^{m+1}_{p(m, k)+i-1, p(m, k)+i}} \times \frac{1}{\Delta^{m+1}_{p(m, k), p(m, k)+i}}\right)^{\frac{1}{2}}, &\quad\text{if } t\in \left[t_{p(m,k)+i-1}^{m+1},t_{p(m,k)+i}^{m+1}\right)
			\end{cases}.
		\end{equation} 
	\end{definition}
	
	We note that the function values of $\psi^{\pi}_{m, k, i}$ are chosen to satisfy $\int \psi^{\pi}_{m, k, i}(t)dt = 0$ and $\int (\psi^{\pi}_{m, k, i}(t))^2 dt = 1$ so that the collection $\{\psi^{\pi}_{m, k, i}\}$ is an orthonormal basis in $L^2([0, T])$. The Schauder functions $e^{\pi}_{m,k,i} : [0,T] \rightarrow \mathbb{R}$ are obtained by integrating the generalized Haar basis:
	\begin{equation*}
		e^{\pi}_{m,k,i}(t) := \int_0^t \psi^{\pi}_{m,k,i}(s)ds =\left(\int_{t^{m+1}_{p(m,k)}}^{t\wedge t^{m+1}_{p(m,k)+i}} \psi^{\pi}_{m,k,i}(s)ds\right) \mathbbm{1}_{[t^m_k,t^{m+1}_{p(m,k)+i}]}(t). 
	\end{equation*}
	
	To further simplify the notations in what follows, we introduce
	\begin{align*}
		&t^{m,k,i}_1 := t^{m+1}_{p(m, k)}, \qquad t^{m,k,i}_2 := t^{m+1}_{p(m, k)+i-1}, \qquad t^{m,k,i}_3 := t^{m+1}_{p(m, k)+i},
		\\ 
		\Delta^{m,k,i}_1 &:= \Delta^{m+1}_{p(m, k), p(m, k)+i-1} = t^{m,k,i}_2 - t^{m,k,i}_1, \qquad
		\Delta^{m,k,i}_2 := \Delta^{m+1}_{p(m, k)+i-1, p(m, k)+i} = t^{m,k,i}_3 - t^{m,k,i}_2.
	\end{align*}
	
	\begin{definition}  [Schauder basis] \label{Def : Schauder functions}
		For every index $m,k,i$ of Definition~\ref{Def : Haar}, the following function $e^\pi_{m,k,i}$ is called \textit{generalized Schauder function} associated with $\pi = (\pi^n)_{n \ge 0}$:
		\begin{eqnarray} \label{defn.e}
			e^{\pi}_{m,k,i}(t) =
			\begin{cases}
				\qquad \qquad \qquad \qquad 0, & ~~ \text{if } t\notin [t_1^{m, k, i}, t_3^{m, k, i}) 
				\\
				\left( \frac{\Delta^{m, k, i}_2}{\Delta^{m, k, i}_{1}} \times \frac{1}{\Delta^{m, k, i}_{1} + \Delta^{m, k, i}_{2}} \right)^{\frac{1}{2}}\times(t-t^{m,k,i}_1), & ~~ \text{if } t\in [t_1^{m,k,i}, t_2^{m,k,i})
				\\
				\left( \frac{\Delta^{m,k,i}_{1}}{\Delta^{m,k,i}_{2}} \times \frac{1}{\Delta^{m,k,i}_{1}+ \Delta^{m,k,i}_2}\right)^{\frac{1}{2}} \times(t^{m,k,i}_3-t), & ~~\text{if } t \in [t_2^{m,k,i},t_3^{m,k,i})
			\end{cases}.
		\end{eqnarray}
	\end{definition}
	
	Note that generalized Schauder functions are continuous, triangle-shaped (and not differentiable) functions. The following result shows that any continuous function defined on $[0, T]$ admits a unique Schauder representation along a given partition sequence $\pi$.
	
	\begin{proposition} [Theorem 3.8 of \cite{das2021}] \label{prop:coeff_hat_func}
		Let $\pi$ be a finitely refining partition sequence of $[0,T]$. Then, every continuous function $x :[0,T] \rightarrow \mathbb{R}$ has a unique Schauder representation along $\pi$:
		\begin{equation}    \label{Schauder representation}
			x(t) = x(0) + \big(x(T)-x(0)\big)t+ \sum_{m=0}^{\infty} \sum_{k=0}^{N(\pi^m)-1} \sum_{i = 1}^{p(m, k+1)-p(m, k)} \theta^{x, \pi}_{m,k,i} e^{\pi}_{m,k,i}(t), \qquad \forall \, t \in [0, T],
		\end{equation}
		with a closed-form representation of the Schauder coefficient
		\begin{equation}  \label{eq.theta.coeff}
			\theta^{x, \pi}_{m,k,i} = \frac{\big(x(t^{m,k,i}_{2})-x(t^{m,k,i}_{1})\big)(t^{m,k,i}_{3}-t^{m,k,i}_{2})-\big(x(t^{m,k,i}_{3})- x(t^{m,k,i}_{2})\big)(t^{m,k,i}_{2}-t^{m,k,i}_1)}{\sqrt{(t^{m,k,i}_2-t^{m,k,i}_1)(t^{m,k,i}_3-t^{m,k,i}_2)(t^{m,k,i}_3-t^{m,k,i}_1)}}.
		\end{equation}
	\end{proposition}
	
\begin{remark}	\label{rem: re-indexing}
A family of Schauder functions $\{e^{\pi}_{m,k,i}\}_{m,k,i}$ in Definition \ref{Def : Schauder functions} can be reordered as $\{e^{\pi}_{m,k}\}_{m,k}$, such that for each $m \ge 0$ the values of $k$ run from $0$ to $N(\pi^{m+1})-N(\pi^m)-1$ after reordering. We shall frequently use this reordering to simplify the notation and denote the index set 
\begin{equation}    \label{def : index set}
I_m := \{ 0, 1, \cdots, N(\pi^{m+1})-N(\pi^{m})-1\} 
\end{equation}
for each $m$. The corresponding Schauder coefficients $\{\theta^{x, \pi}_{m,k,i}\}_{m,k,i}$ in Proposition \eqref{prop:coeff_hat_func} can be reordered as $\{\theta^{x, \pi}_{m,k}\}_{m,k}$ for $k \in I_m$ and $m \ge 0$ in the same manner.
\end{remark}

\begin{remark}
Consider the following approximation $x_n$ up to level $n$ of $x$ in \eqref{Schauder representation}
\begin{equation*}
    x_n(t):= x(0) + \big(x(T)-x(0)\big)t+ \sum_{m=0}^{n-1} \sum_{k \in I_m} \theta^{x, \pi}_{m,k,i} e^{\pi}_{m,k,i}(t), \qquad \forall \, t \in [0, T].
\end{equation*}
Then, it is straightforward to show that $x_n$ is just the piecewise linear approximation of $x$ along the partition points of $\pi^n$, i.e.,
\begin{align*}
    x_n(s) =
    \begin{cases}
        ~~~~~~~~~ x(s), \quad &\text{if } s \in \pi^n,
        \\
        x(t^n_i) + \frac{x(t^n_{i+1})- x(t^n_{i})}{t^n_{i+1}-t^n_i}\times(s-t^n_i), \quad &\text{if } s \in (t^n_i,t^n_{i+1}) \text{ for some } t^n_i \in \pi^n.
    \end{cases}
\end{align*}
\end{remark}

\bigskip

\section{Characterization of variation index}  \label{sec: main result}
In this section, we characterize the variation index $p^{\pi}(x)$ of $x \in C^0([0, T])$ along $\pi \in \Pi([0, T])$, in terms of the Schauder coefficients $\{\theta^{x, \pi}_{m,k}\}_{m,k}$ introduced in Section~\ref{subsec: Haar and Schauder}. This characterization provides an equivalent condition for a continuous function $x$ to belong to the Banach space $\mathcal{X}^p_{\pi}$. We note that the Schauder basis expansions have already been used to characterize different function spaces that provide regularity properties of functions, such as H\"older space \cite{fake_fBM,Ciesielski:isomorphism} and Besov space \cite{rosenbaum2009, ciesielski1993}. We recall the definition \eqref{def: discrete p-th variation} of the $p$-th variation, as well as Definitions \ref{def.finite.refining}-\ref{def.conv.ref}.
	
    \begin{remark}  \label{rem : normalizing x(0)}
		Any $x \in C^0([0, T])$ can be translated to $\bar{x} \in C^0([0, T])$ with $\bar x(0) = \bar x(T) = 0$, by adding a linear function. For any $p > 1$, the $p$-th variation of a linear function $y$ along any element $\pi = (\pi^n)_{n \ge 0}$ of $\Pi([0, T])$ is zero, i.e., $\limsup_{n \to \infty} [y]^{(p)}_{\pi^n} = 0$. Moreover, the subadditive property of the norm $\Vert \cdot \Vert^{(p)}_{\pi}$ in Definition \ref{Def : p norm and Xp space} implies $\Vert \bar x\Vert_\pi^{(p)} < \infty$ if and only if $\Vert x \Vert_\pi^{(p)} < \infty$. Since we are only interested in the conditions regarding the finiteness of $\Vert x \Vert^{(p)}_{\pi}$-norm (or $\limsup_{n \to \infty} [x]^{(p)}_{\pi^n}$), we shall assume without loss of generality $x(0) = x(T) = 0$ in what follows. Then, the Schauder representation \eqref{Schauder representation} of any $x \in C^0([0, T])$ becomes simpler:
		\begin{equation}	\label{eq : simpler Schauder representation}
		  x(t) = \sum_{m=0}^{\infty} \sum_{k=0}^{N(\pi^m)-1} \sum_{i = 1}^{p(m, k+1)-p(m, k)} \theta^{x, \pi}_{m,k,i} e^{\pi}_{m,k,i}(t), \qquad \forall \, t \in [0, T].
		\end{equation}
		The above triple sum can be expressed as a double sum after re-indexing as in Remark \ref{rem: re-indexing}.
	\end{remark}

	\medskip
	
\subsection{Results} 	\label{subsec: variation index results}	
We provide Proposition \ref{prop : p-th var general} and Theorem \ref{thm: p-th var general2} below, and their proofs are given in the next subsection. 
	
\begin{proposition} \label{prop : p-th var general} 
		For any $p > 1$, $x \in C^0([0, T])$, and a balanced, complete refining partition sequence $\pi = (\pi^n)_{n \ge 0}$ of $[0, T]$, we denote 
		\begin{equation}   \label{def: eta_n general}
			\eta^{\pi, (p)}_n := |\pi^n|^{p-1} \Bigg( \sum_{m=0}^{n-1} \vert \pi^{m} \vert^{\frac{1}{p}-\frac{1}{2}} \bigg( \sum_{k \in I_m} |\theta^{x, \pi}_{m, k}|^p \bigg)^\frac{1}{p} \Bigg)^{p}.
		\end{equation}
		Then, we have
		\begin{equation}    \label{ineq : p-th var, eta limsup}
			\limsup_{n \to \infty} \, [x]_{\pi^n}^{(p)}(T) < \infty \quad \text{if and only if} \quad \limsup_{n \to \infty} \, \eta^{\pi, (p)}_n < \infty.
		\end{equation}
	\end{proposition}

For any balanced, complete refining partition sequence $\pi$, Proposition \ref{prop : p-th var general} immediately provides the sufficient and necessary condition for $x \in C^0([0, T])$ to belong to the Banach space $\mathcal{X}^p_{\pi}$ in \eqref{def : x p pi space}, in terms of its Schauder coefficients through the real sequence $(\eta^{\pi, (p)}_n)_{n \ge 0}$:
	\begin{equation*}
		x \in \mathcal{X}^p_{\pi} \quad \Longleftrightarrow \quad \limsup_{n \to \infty} \, \eta^{\pi, (p)}_n < \infty.
	\end{equation*}
	Moreover, it also yields the equivalent formulation of the variation index in \eqref{def. alter variation index}:
	\begin{equation}    \label{def. alter variation index2 general}
		p^{\pi}(x) = \inf \big\{ p > 1 : \limsup_{n \to \infty} \, \eta^{\pi, (p)}_n < \infty \big\}.
	\end{equation}

	Thus, the $(\limsup)$-finiteness of the sequence $(\eta^{\pi, (p)}_n)_{n \ge 0}$ can provide useful path property of $x$ along any balanced, complete refining partition sequences, and each term $\eta^{\pi, (p)}_n$ contains the Schauder coefficients of $x$ up to level $n-1$, namely $\{\theta^{x, \pi}_{m, k}\}_{m = 0, \cdots, n-1, \, k \in I_m}$. However, with nominal additional conditions on the partition sequence, we have a simpler condition involving Schauder coefficients.	
	
	\begin{theorem} \label{thm: p-th var general2}
		For any $p > 1$, $x \in C^0([0, T])$, and a balanced, convergent refining partition sequence $\pi = (\pi^n)_{n \ge 0}$ of $[0, T]$, we denote 
		\begin{equation}    \label{def : xi general}
			\xi_n^{\pi,(p)} = |\pi^n|^{\frac{p}{2}} \bigg( \sum_{k \in I_n} |\theta^{x, \pi}_{n, k}|^p \bigg), \qquad \forall \, n \ge 0.
		\end{equation}
		Then, we have
		\begin{equation}	\label{ineq : p-th var, xi limsup}
			\limsup_{n \to \infty} \, [x]_{\pi^n}^{(p)}(T) < \infty \quad \text{if and only if} \quad \limsup_{n \to \infty} \, \xi^{\pi, (p)}_n < \infty.
		\end{equation}
		Thus, we also have
		\begin{equation*}
			x \in \mathcal{X}^p_{\pi} \quad \text{if and only if} \quad \limsup_{n \to \infty} \, \xi^{\pi, (p)}_n < \infty.
		\end{equation*}
	\end{theorem}
	
	In the definition \eqref{def : xi general}, the quantity $\xi^{\pi, (p)}_n$ only contains the Schauder coefficients $\{\theta^{x, \pi}_{n, k}\}_{k \in I_n}$ of $x$ that belong to the $n$-th level, for each $n \in \N$.
	Theorem \ref{thm: p-th var general2} also provides a similar equivalent formulation of the variation index in \eqref{def. alter variation index}.
	
\begin{corollary}  \label{cor : new def variation index general2}
		Let $\pi$ be a balanced, convergent refining partition sequence. Then, we have
		\begin{equation}    \label{def. alter variation index3 general}
			p^{\pi}(x) = \inf \big\{ p > 1 : \limsup_{n \to \infty} \, \xi^{\pi, (p)}_n < \infty \big\}.
		\end{equation}
\end{corollary}

\begin{remark}
    In all of the previous results, we considered the (generalized) $p$-th variation up to the terminal time $T$. However, we can derive similar results for any points $t \in [0,T]$. For $x \in C^0([0, T])$, let us recall the definition \eqref{def: p-th variation} of $[x]^{(p)}_{\pi^n}(t)$ such that the mapping $t \mapsto \limsup_{n \to \infty} \, [x]^{(p)}_{\pi^n}(t)$ is nondecreasing. We also introduce the notations
    \begin{align}   
        \eta^{\pi, (p)}_n (t)&:= |\pi^n|^{p-1} \Bigg( \sum_{m=0}^{n-1} \vert \pi^{m} \vert^{\frac{1}{p}-\frac{1}{2}} \bigg( \sum_{\substack{k \in I_m\\ supp(e^{\pi}_{m,k}) \subset [0,t]}} |\theta^{x, \pi}_{m, k}|^p \bigg)^\frac{1}{p} \Bigg)^{p},     \label{def: eta_n general (t)}
        \\
        \xi_n^{\pi,(p)}(t) &:= |\pi^n|^{\frac{p}{2}} \bigg( \sum_{\substack{k \in I_n\\ supp(e^{\pi}_{n,k}) \subset [0,t]}} |\theta^{x, \pi}_{n, k}|^p \bigg).            \label{def : xi general (t)}
    \end{align}
    Then, the results \eqref{ineq : p-th var, eta limsup} and \eqref{ineq : p-th var, xi limsup} can be replaced by
    \begin{align}
        \limsup_{n \to \infty} \, [x]_{\pi^n}^{(p)}(t) < \infty \quad &\text{if and only if} \quad \limsup_{n \to \infty} \, \eta^{\pi, (p)}_n (t) < \infty, \quad \text{and}        \label{ineq : p-th var, eta limsup (t)}
        \\
        \limsup_{n \to \infty} \, [x]_{\pi^n}^{(p)}(t) < \infty \quad &\text{if and only if} \quad \limsup_{n \to \infty} \, \xi^{\pi, (p)}_n(t) < \infty, \quad \text{for every } t \in \cup_{n \in \N} \pi^n.          \label{ineq : p-th var, xi limsup (t)}
    \end{align}
    To show \eqref{ineq : p-th var, eta limsup (t)} and \eqref{ineq : p-th var, xi limsup (t)}, we first define a `stopped function' $x_t(s) := x(t \wedge s)$ for $s \in [0, T]$. Furthermore, we define \begin{equation*}
        \widetilde\theta^{x, \pi}_{m,k} :=
        \begin{cases}
            \theta^{x, \pi}_{m,k}, &\text{ if supp}(e^\pi_{m,k}) \subset [0,t],
            \\
            ~0, &\text{ otherwise,}
        \end{cases}
    \end{equation*} 
    and
    \begin{equation*}
        \widetilde x(t):=\sum_{m=0}^\infty \sum_{k \in I_m} \widetilde\theta^{x, \pi}_{m,k}e^\pi_{m,k}(t).
    \end{equation*}
    For $t \in \cup_{n \in \N} \pi^n =: P$, the two functions $x_t$ and $\widetilde x$ differ only by a finite sum of piecewise linear functions, say $y$, which hence satisfies $[y]^{(p)}_{\pi} \equiv 0$ for every $p>1$. Lemma \ref{lem: adding zero p-th variation} therefore yields that $\limsup_{n \to \infty} \, [\widetilde x]_{\pi^n}^{(p)}(T) =\limsup_{n \to \infty} \, [x_t]^{(p)}_{\pi^n}(T) = \limsup_{n \to \infty} \, [x]_{\pi^n}^{(p)} (t)$. Now applying Proposition \ref{prop : p-th var general} and Theorem \ref{thm: p-th var general2} to $\widetilde{x}$ with the quantities \eqref{def: eta_n general (t)} and \eqref{def : xi general (t)}, proves \eqref{ineq : p-th var, eta limsup (t)} and \eqref{ineq : p-th var, xi limsup (t)}. 
    
    For $t \notin P$, we can choose a point $s \in P$ which is sufficiently close and bigger than $t$, and check the finiteness of $\limsup_{n \to \infty} \eta^{\pi, (p)}_{n}(s)$, or $\limsup_{n \to \infty} \xi^{\pi, (p)}_{n}(s)$, to conclude the finiteness $\limsup_{n \to \infty} [x]^{(p)}_{\pi^n}(t) \le \limsup_{n \to \infty} [x]^{(p)}_{\pi^n}(s) < \infty$.
\end{remark}

	\medskip
	
	\subsection{Proofs}
	
	Before proving Proposition \ref{prop : p-th var general} and Theorem \ref{thm: p-th var general2}, we first introduce some preliminary lemmata.
	\begin{lemma}   \label{lem: a_n/b_n}
		Let $(a_n)_{n \in \N}$ and $(b_n)_{n \in \N}$ be real sequences such that $b_n > 0$, $\frac{b_{n+1}}{b_n} =: \beta_n > 1$ for every $n \in \N$, and the limit $\lim_{n \to \infty} \beta_n = \beta > 1$ exists. Then, we have the inequality
		\begin{equation}	\label{ineq : limsup}
			\limsup_{n \rightarrow \infty} \bigg( \frac{a_{n+1} - a_n}{b_{n+1} - b_n} \bigg)  \le \frac{\beta}{\beta-1} \limsup_{n \rightarrow \infty} \bigg( \frac{a_{n+1}}{b_{n+1}} \bigg) - \frac{1}{\beta-1} \liminf_{n \rightarrow \infty} \bigg( \frac{a_n}{b_n} \bigg).
		\end{equation}
	\end{lemma}
	
	\begin{proof}[Proof of Lemma \ref{lem: a_n/b_n}]
		Taking $\limsup$ to the both sides of the following identity
		\begin{equation}    \label{eq : a_n b_n identity}
			\frac{a_{n+1} - a_n}{b_{n+1} - b_n} = \frac{1}{\frac{b_{n+1}}{b_n} - 1}\bigg( \frac{a_{n+1}}{b_{n+1}} \times \frac{b_{n+1}}{b_n} - \frac{a_n}{b_n}\bigg)
			= \frac{1}{\beta_n-1} \bigg( \beta_n\frac{a_{n+1}}{b_{n+1}} - \frac{a_n}{b_n} \bigg)
		\end{equation}
		with the following properties for any real sequences $(x_n)_{n \in \N}, (y_n)_{n \in \N}$ proves the result:
		\begin{align}
			&\limsup_{n \to \infty} \, (x_n+y_n) \le \limsup_{n \to \infty} x_n + \limsup_{n \to \infty} y_n, \qquad
			\limsup_{n \to \infty} \, (-x_n) = -\liminf_{n \to \infty} x_n,    \label{eq : -limsup}
			\\
			&\limsup_{n \to \infty} \, (x_n y_n) = \big( \lim_{n \to \infty} x_n \big) \big(\limsup_{n \to \infty} y_n\big), \quad \text{provided that } \lim_{n \to \infty} x_n  \text{ exists and is positive.} \nonumber
		\end{align}    
	\end{proof}
	
	\begin{lemma}   \label{lem : liminf and limsup}
		Let $(a_n)_{n \in \N}$ and $(b_n)_{n \in \N}$ be real sequences such that $(b_n)_{n \in \N}$ is strictly increasing and $\lim_{n \to \infty} b_n = \infty$. Then, we have the following inequalities
		\begin{equation}    \label{ineq : limsup and liminf}
			\liminf_{n \rightarrow \infty} \bigg( \frac{a_{n+1} - a_n}{b_{n+1} - b_n} \bigg)
			\le
			\liminf_{n \rightarrow \infty} \bigg( \frac{a_n}{b_n} \bigg)
			\le 
			\limsup_{n \rightarrow \infty} \bigg( \frac{a_n}{b_n} \bigg)
			\le
			\limsup_{n \rightarrow \infty} \bigg( \frac{a_{n+1} - a_n}{b_{n+1} - b_n} \bigg).
		\end{equation}
	\end{lemma}
	
	\begin{proof}[Proof of Lemma \ref{lem : liminf and limsup}]
		The middle inequality is obvious. We shall show the last inequality; the first inequality then follows from \eqref{eq : -limsup}. If the right-most term of \eqref{ineq : limsup and liminf} diverges to infinity, there is nothing to show. Thus, we assume
		\begin{equation*}
			\limsup_{n \rightarrow \infty} \bigg( \frac{a_{n+1} - a_n}{b_{n+1} - b_n} \bigg) = L < \infty.
		\end{equation*}
		For any $r > L$, there exists $N \in \N$ such that 
		\begin{equation*}
			\frac{a_{n+1} - a_n}{b_{n+1} - b_n} < r, \qquad \text{or} \qquad a_{n+1} - a_n < r (b_{n+1}-b_n),
		\end{equation*}
		holds for every $n > N$. Fix an arbitrary integer $m$ greater than $N$, and sum up the last inequalities for $n = N, \cdots, m-1$ to obtain
		\begin{equation*}
			a_m - a_N = \sum_{n=N}^{m-1} (a_{n+1} - a_n) < r\sum_{n=N}^{m-1} (b_{n+1} - b_n) = r (b_{m}-b_N), \quad \text{thus} \quad \frac{a_m-a_N}{b_m} < r - r\frac{b_N}{b_m}.
		\end{equation*}
		Sending $m$ to infinity and using the fact $\lim_{m \to \infty} b_m = \infty$ yields the inequality
		\begin{equation*}
			\limsup_{m \rightarrow \infty} \bigg( \frac{a_m}{b_m} \bigg) < r.
		\end{equation*}
		Since this should hold for any $r > L$, we conclude that the last inequality of \eqref{ineq : limsup and liminf} holds.
	\end{proof}
	
	\begin{lemma}	\label{lem : Silverman-Toeplitz weak}
		Let $A = (a_{n, m})_{n \ge 0, m \ge 0}$ be an infinite-dimensional matrix satisfying the following properties:
		\begin{enumerate} [(i)]
			\item $\lim_{n \rightarrow \infty} a_{n, m} = 0$ for every $m \ge 0$;
			\item $\lim_{n \rightarrow \infty} \sum_{m=0}^{\infty} a_{n, m} = 1$;
			\item $\sup_{n \ge 0} \sum_{m=0}^{\infty} |a_{n, m}| < \infty$.
		\end{enumerate}
		Then, for any real sequence $(s_n)_{n \ge 0}$ with nonnegative terms, i.e., $s_n \ge 0$ for all $n \ge 0$, we have 
		\begin{equation}\label{eq: silverman gen}
			\limsup_{n \to \infty} \sum_{m=0}^{\infty} a_{n, m} s_m \le \limsup_{n \to \infty} s_n.
		\end{equation} 
	\end{lemma}
	
	\begin{remark}
		We note that Lemma \ref{lem : Silverman-Toeplitz weak} is inspired by the Silverman-Toeplitz Theorem (see, e.g., \citep{boos2000}), which states that the real sequence $(s_n)_{n \ge 0}$ converges to $s$, if and only if 
		\begin{equation}\label{eq: silverman}
			\lim_{n \to \infty} \Big( \sum_{m=0}^{n} a_{n, m} s_m \Big) = s,
		\end{equation}
		for $A = (a_{n, m})_{n \ge 0, m \ge 0}$ satisfying the conditions of Lemma \ref{lem : Silverman-Toeplitz weak}.
	\end{remark}

	\begin{proof}[Proof of Lemma \ref{lem : Silverman-Toeplitz weak}]
		If $\limsup_{n \to \infty} s_n = \infty$, then there is nothing to prove; thus, we assume $\limsup_{n \to \infty} s_n =: s < \infty$. This implies that there exists $K < \infty$ such that $s_n \le K$ for all $n \ge 0$. 
		We denote $L := \sup_{n \ge 0} \sum_{m=0}^{\infty} |a_{n, m}| < \infty$ in condition (iii), and fix an arbitrary $\epsilon > 0$. Then, there exists $M_1 \in \mathbb{N}$ such that
		\begin{equation}    \label{sm bound}
			s_m \le s + \frac{\epsilon}{4L}, \qquad \text{for every } m > M_1.
		\end{equation}
		Condition (i) implies that there exist constants $N_0, N_1, \cdots, N_{M_1}$ such that
		\begin{equation*}
			|a_{n, m}| \le \frac{\epsilon}{4(M_1+1)(K+1)}, \qquad \text{for every } 0 \le m \le M_1 \text{ and } n > N_{m}.
		\end{equation*}
		Set $\tilde{N} := \max \{N_0, N_1, \cdots, N_{M_1}\}$, then
		\begin{equation*}
			\sum_{m=0}^{M_1} a_{n, m}s_m \le \sum_{m=0}^{M_1} |a_{n, m}s_m| \le \sum_{m=0}^{M_1} \frac{s_m \epsilon}{4(M_1+1)(K+1)} < \frac{\epsilon}{4}, \qquad \text{for every } n > \tilde{N}.
		\end{equation*}
		On the other hand, we have from \eqref{sm bound}
		\begin{equation*}
			\sum_{m=M_1+1}^{\infty} a_{n, m}s_m 
			\le s \sum_{m=M_1+1}^{\infty} |a_{n, m}| + \frac{\epsilon}{4L} \sum_{m=M_1+1}^{\infty} |a_{n, m}|
			\le s \sum_{m=M_1+1}^{\infty} |a_{n, m}| + \frac{\epsilon}{4}.
		\end{equation*}
		Combining the last two inequalities,
		\begin{equation}    \label{two sums}
			\sum_{m=0}^{\infty} a_{n, m}s_m
			= \sum_{m=0}^{M_1} a_{n, m}s_m + \sum_{m=M_1+1}^{\infty} a_{n, m}s_m
			\le s\sum_{m=M_1+1}^{\infty} |a_{n, m}| + \frac{\epsilon}{2} \quad \text{for every } n > \tilde{N}.
		\end{equation}
		We now claim that $(\sum_{m=0}^{\infty} a_{n, m}s_m)_{n \ge 0}$ is an absolutely convergence sequence
		\begin{equation*}
			\sum_{m = 0}^{\infty} |a_{n, m}s_m| \le K \sum_{m=0}^{\infty} |a_{n, m}| \le K L < \infty,
		\end{equation*}
		thanks to condition (iii). Therefore, taking the limit as $n \to \infty$ in \eqref{two sums}, together with condition (ii), we conclude
		\begin{equation*}
			\limsup_{n \to \infty} \sum_{m=0}^{\infty} a_{n, m}s_m \le s + \frac{\epsilon}{2}.
		\end{equation*}
		Since $\epsilon$ is chosen arbitrarily, this proves the result.
	\end{proof}
	
	We are now ready to prove Proposition \ref{prop : p-th var general} and Theorem \ref{thm: p-th var general2}.
	
	\begin{proof} [Proof of Proposition \ref{prop : p-th var general}]
		Using the Schauder representation \eqref{eq : simpler Schauder representation}, we expand the $p$-th variation of $x$ along $\pi^n$ for each $n \in \mathbb{N}$,
		\begin{align}
			[x]_{\pi^n}^{(p)}(T) &= \sum_{\ell=0}^{N(\pi^n)-1} \Big| x(t^n_{\ell+1})-x(t^n_{\ell}) \Big|^p 	\label{eq: p-th variation expression general}    
			\\
			&= \sum_{\ell=0}^{N(\pi^n)-1} \bigg|\sum_{m=0}^{n-1} \sum_{k=0}^{N(\pi^m)-1} \sum_{i = 1}^{p(m, k+1)-p(m, k)} \theta^{x, \pi}_{m,k,i} \Big( e_{m,k,i}(t^n_{\ell+1})-e_{m,k,i}(t^n_{\ell}) \Big) \bigg|^p.		\nonumber
		\end{align}
		Since $\pi$ is finitely refining, for each fixed pair $(m, \ell)$ with $m < n$ and $\ell < N(\pi^n)$, the cardinality of the set $I(m, \ell) := \{(k, i) : e_{m, k,i}(t^n_{\ell+1}) - e_{m, k,i}(t^n_{\ell}) \neq 0\}$ has an upper bound $M$. Also, in Definition \ref{Def : Schauder functions}, we note that 
		\begin{equation*}
			\underline{\pi^{m+1}} \le \Delta^{m,k,i}_1 \le M|\pi^{m+1}|, \qquad
			\underline{\pi^{m+1}} \le \Delta^{m,k,i}_2 \le |\pi^{m+1}|,
		\end{equation*}
		as $\Delta^{m,k,i}_1$ is a length of an interval containing at most $M$ many consecutive intervals of $\pi^{m+1}$, whereas $\Delta^{m,k,i}_2$ is a length of a single interval of $\pi^{m+1}$. From the balanced and complete refining property, we have
		\begin{align*}
			\Big\vert e_{m, k, i}(t^n_{\ell+1}) - e_{m, k, i}(t^n_{\ell}) \Big\vert
			&\le \frac{1}{\sqrt{\Delta^{m,k,i}_1 + \Delta^{m,k,i}_2}}\Bigg( \max \bigg(\sqrt{\frac{\Delta^{m,k,i}_2}{\Delta^{m,k,i}_1}}, \sqrt{\frac{\Delta^{m,k,i}_1}{\Delta^{m,k,i}_2}} \bigg) \Bigg) \vert \pi^n \vert
			\\
			&\le \frac{1}{\sqrt{\underline{\pi^{m+1}}}} \sqrt{\frac{M|\pi^{m+1}|}{\underline{\pi^{m+1}}}} \vert \pi^n \vert
			\le \frac{\sqrt{cM}}{\sqrt{\underline{\pi^{m+1}}}}  \vert \pi^n \vert
			\le \frac{c\sqrt{M}}{\sqrt{\vert \pi^{m+1} \vert}} \vert \pi^n \vert = \frac{c\sqrt{bM} \vert \pi^n \vert}{\sqrt{\vert \pi^{m} \vert}}.
		\end{align*}
		Thus, we have from \eqref{eq: p-th variation expression general}
		\begin{align*}
			[x]_{\pi^n}^{(p)}(T) &\le \sum_{\ell=0}^{N(\pi^n)-1} \Bigg\vert \sum_{m=0}^{n-1} M \Big( \max_{(k, i) \in I(m, \ell)} |\theta^{x, \pi}_{m, k, i}| \Big) \frac{c\sqrt{bM} \vert \pi^n \vert}{\sqrt{\vert \pi^m \vert}} \Bigg\vert^p
			\\
			&= \Big( Mc\sqrt{bM}|\pi^n| \Big)^{p} \sum_{\ell=0}^{N(\pi^n)-1}  \Bigg| \sum_{m=0}^{n-1} \Big( \max_{(k, i) \in I(m, \ell)} |\theta^{x, \pi}_{m, k, i}| \Big) \vert \pi^m \vert^{-\frac{1}{2}} \Bigg|^p =: Q_n.
		\end{align*}
		We now set $\epsilon := p - \lfloor p \rfloor$ and expand the $\lfloor p \rfloor$-th power to obtain
		\begin{align*}
			&\frac{Q_n}{\big(Mc\sqrt{bM}|\pi^n|\big)^{p}}
			= \sum_{\ell=0}^{N(\pi^n)-1} \left| \sum_{m=0}^{n-1} \Big( \max_{(k, i) \in I(m, \ell)} |\theta^{x, \pi}_{m, k, i}| \Big) \vert \pi^m \vert^{-\frac{1}{2}} \right|^{\lfloor p \rfloor} \left| \sum_{m=0}^{n-1} \Big( \max_{(k, i) \in I(m, \ell)} |\theta^{x, \pi}_{m, k, i}| \Big) \vert \pi^m \vert^{-\frac{1}{2}} \right|^{\epsilon}
			\\
			&= \sum_{\ell=0}^{N(\pi^n)-1} \sum_{0 \le m_1, \cdots, m_{\lfloor p \rfloor} \le n-1} \left( \prod_{j=1}^{\lfloor p \rfloor} \Big( \max_{(k, i) \in I(m_j, \ell)} |\theta^{x, \pi}_{m_j, k, i}| \Big) \vert \pi^{m_j} \vert^{-\frac{1}{2}} \right) \left| \sum_{m=0}^{n-1} \Big( \max_{(k, i) \in I(m, \ell)} |\theta^{x, \pi}_{m, k, i}| \Big) \vert \pi^m \vert^{-\frac{1}{2}} \right|^{\epsilon}
			\\
			&= \sum_{0 \le m_1 ,\cdots, m_{\lfloor p \rfloor} \le n-1} \bigg( \prod_{j=1}^{\lfloor p \rfloor} \vert \pi^{m_j} \vert^{-\frac{1}{2}} \bigg) \sum_{\ell=0}^{N(\pi^n)-1} \bigg( \prod_{j=1}^{\lfloor p \rfloor} \max_{(k, i) \in I(m_j, \ell)} |\theta^{x, \pi}_{m_j, k, i}| \bigg) \left| \sum_{m=0}^{n-1} \Big( \max_{(k, i) \in I(m, \ell)} |\theta^{x, \pi}_{m, k, i}| \Big) \vert \pi^m \vert^{-\frac{1}{2}} \right|^{\epsilon}
			\\
			& \le \sum_{0 \le m_1 ,\cdots, m_{\lfloor p \rfloor} \le n-1} \bigg( \prod_{j=1}^{\lfloor p \rfloor} \vert \pi^{m_j} \vert^{-\frac{1}{2}} \bigg) 
			\\
			& \qquad \qquad \qquad \times \prod_{j=1}^{\lfloor p \rfloor} \bigg( \sum_{\ell=0}^{N(\pi^n)-1} \max_{(k, i) \in I(m_j, \ell)} |\theta^{x, \pi}_{m_j, k, i}|^p \bigg)^\frac{1}{p} \bigg( \sum_{\ell=0}^{N(\pi^n)-1} \left| \sum_{m=0}^{n-1}  \Big( \max_{(k, i) \in I(m, \ell)} |\theta^{x, \pi}_{m, k, i}| \Big) \vert \pi^m \vert^{-\frac{1}{2}} \right|^{\epsilon \cdot \frac{p}{\epsilon}} \bigg)^{\frac{\epsilon}{p}}
			\\
			& = \sum_{0 \le m_1 ,\cdots, m_{\lfloor p \rfloor} \le n-1} \bigg( \prod_{j=1}^{\lfloor p \rfloor} \vert \pi^{m_j} \vert^{-\frac{1}{2}} \bigg) \prod_{j=1}^{\lfloor p \rfloor} \bigg( \sum_{\ell=0}^{N(\pi^n)-1} \max_{(k, i) \in I(m_j, \ell)} |\theta^{x, \pi}_{m_j, k, i}|^p \bigg)^\frac{1}{p} \bigg( \frac{Q_n}{\big(Mc\sqrt{bM} |\pi^n|\big)^{p}} \bigg)^{\frac{\epsilon}{p}}.
		\end{align*}
		Here, the inequality follows from generalized H\"older inequality with $\frac{1}{p} \times \lfloor p \rfloor + \frac{\epsilon}{p} = 1$. We further derive
		\begin{align*}
			(Q_n)^{1-\frac{\epsilon}{p}}
			& \le \Big(Mc\sqrt{bM}|\pi^n|\Big)^{\lfloor p \rfloor} \sum_{0 \le m_1 \cdots m_{\lfloor p \rfloor} \le n-1} \bigg( \prod_{j=1}^{\lfloor p \rfloor} \vert \pi^{m_j} \vert^{-\frac{1}{2}} \bigg) \prod_{j=1}^{\lfloor p \rfloor} \bigg( \sum_{\ell=0}^{N(\pi^n)-1} \max_{(k, i) \in I(m_j, \ell)} |\theta^{x, \pi}_{m_j, k, i}|^p \bigg)^\frac{1}{p}
			\\
			& \le \Big(Mc\sqrt{bM}|\pi^n|\Big)^{\lfloor p \rfloor} \sum_{0 \le m_1 \cdots m_{\lfloor p \rfloor} \le n-1} \bigg( \prod_{j=1}^{\lfloor p \rfloor} \vert \pi^{m_j} \vert^{-\frac{1}{2}} \bigg) \prod_{j=1}^{\lfloor p \rfloor} \bigg( \frac{c|\pi^{m_j}|}{|\pi^n|} \sum_{k, i} |\theta^{x, \pi}_{m_j, k, i}|^p \bigg)^\frac{1}{p}
			\\
			& = \Big(Mc\sqrt{bM}|\pi^n|\Big)^{\lfloor p \rfloor} \Bigg( \sum_{m=0}^{n-1} \vert \pi^{m} \vert^{-\frac{1}{2}} \bigg( \frac{c|\pi^{m}|}{|\pi^n|}  \bigg)^\frac{1}{p} \bigg( \sum_{k, i} |\theta^{x, \pi}_{m, k, i}|^p \bigg)^\frac{1}{p} \Bigg)^{\lfloor p \rfloor}.
		\end{align*}
Here, the second inequality uses the fact that for a fixed $m_j$ there are at most $\frac{|\pi^{m_j}|}{\underline{\pi^n}}$ many partition points of $\pi^n$ sharing the same $\theta^{x, \pi}_{m_j, k, i}$, and this number is bounded by $\frac{c|\pi^{m_j}|}{|\pi^n|}$ due to the balanced condition. 
Therefore, we obtain
\begin{align}
[x]_{\pi^n}^{(p)}(T) &\le Q_n = \big(Q_n^{1-\frac{\epsilon}{p}}\big)^{\frac{p}{\lfloor p \rfloor}}		\label{def: eta general}
\\ & \le \Big(Mc\sqrt{bM}|\pi^n|\Big)^p \Bigg( \sum_{m=0}^{n-1} \vert \pi^{m} \vert^{-\frac{1}{2}} \bigg( \frac{c|\pi^{m}|}{|\pi^n|}\bigg)^\frac{1}{p} \bigg( \sum_{k, i} |\theta^{x, \pi}_{m, k, i}|^p \bigg)^\frac{1}{p} \Bigg)^{p}
= c\Big(Mc\sqrt{bM}\Big)^p \eta^{\pi,(p)}_n,    \nonumber
\end{align}
	from the definition \eqref{def: eta_n general} (after re-indexing $k, i$ into $k$ as in Remark \ref{rem: re-indexing}).
		
On the other hand, using the expression \eqref{eq.theta.coeff} of the Schauder coefficients, we obtain the following bound on the $p$-th power of $\theta^{x, \pi}_{m,k,i}$, thanks to the balanced condition 
\begin{align}
|\theta^{x, \pi}_{m,k,i}|^p \le \bigg(\frac{c}{|\pi^{m+1}|}\bigg)^{\frac{3p}{2}}
\bigg| \big(x(t^{m, k, i}_2)-x(t^{m, k, i}_1)\big) & ( t^{m, k, i}_3 - t^{m, k, i}_2) \label{ineq: p power of theta}
\\
&-\big( x(t^{m, k, i}_3)-x(t^{m, k, i}_2) \big) (t^{m, k, i}_2 -  t^{m, k, i}_1) \bigg|^p.	\nonumber
\end{align}
Here, note that $t^{m, k, i}_2$ and $t^{m, k, i}_3$ are consecutive partition points of $\pi^{m+1}$, but $t^{m, k, i}_1$ and $t^{m, k, i}_2$ may not be. Recalling the notations in \eqref{EqFor_p}, we use the telescoping sum
\begin{equation*}
x(t^{m, k, i}_2)-x(t^{m, k, i}_1) = \sum_{j=1}^{i-1} \Big( x(t^{m+1}_{p(m, k)+j})-x(t^{m+1}_{p(m, k)+j-1}) \Big)
\end{equation*}
with the bound $\max\{|t^{m, k, i}_2 - t^{m, k, i}_1|, |t^{m, k, i}_3 - t^{m, k, i}_2|\} \le M |\pi^{m+1}|$, and apply Jensen's inequality to the right-hand side of \eqref{ineq: p power of theta} to obtain
		\begin{align*}
|\theta^{x, \pi}_{m,k,i}|^p 
& \le \bigg(\frac{c}{|\pi^{m+1}|}\bigg)^{\frac{3p}{2}} (i+1)^{p-1} \Bigg( \sum_{j=1}^{i-1} \Big| \big(x(t^{m+1}_{p(m, k)+j})-x(t^{m+1}_{p(m, k)+j-1}) \big) (t^{m, k, i}_3 -  t^{m, k, i}_2) \Big|^p 
	\\
& \qquad \qquad \qquad \qquad \qquad \qquad \qquad \qquad + \Big| \big(x(t^{m, k, i}_3)-x(t^{m, k, i}_2)\big)(t^{m, k, i}_2 -  t^{m, k, i}_1)\Big|^p \Bigg)
			\\
& \le \frac{M^p c^{\frac{3p}{2}} (i+1)^{p-1}}{|\pi^{m+1}|^{\frac{3p}{2}-p}} \bigg( \sum_{j=1}^{i-1} \big| x(t^{m+1}_{p(m, k)+j})-x(t^{m+1}_{p(m, k)+j-1})\big|^p + \big| x(t^{m, k, i}_3)-x(t^{m, k, i}_2)\big|^p \bigg).
\end{align*}
We note that the quantities inside the last big parenthesis is the $p$-th variation of $x$ along the partition points of $\pi^{m+1}$ that belong to the interval $[t^n_k, t^n_{k+1}]$, and these intervals are disjoint for different values of $k$. We now derive the following inequality
		\begin{equation*}
			\sum_{k=0}^{N(\pi^m)-1} \sum_{i = 1}^{p(m, k+1)-p(m, k)} |\theta^{x, \pi}_{m,k,i}|^p
			\le \frac{M^p c^{\frac{3p}{2}} (M+1)^{p-1}}{|\pi^{m+1}|^{\frac{p}{2}}} M[x]^{(p)}_{\pi^{m+1}}(T) 
			< \frac{c^{\frac{3p}{2}} (M+1)^{2p}}{|\pi^{m+1}|^{\frac{p}{2}}} [x]^{(p)}_{\pi^{m+1}}(T),
		\end{equation*}
		since the largest value $i$ can take is $p(m, k+1)-p(m, k) \le M$ and the first $p$-th power increment $|x(t^{m+1}_{p(m, k)+1})-x(t^{m+1}_{p(m, k)})|^p$ (which has been most repeatedly added) has been added at most $M$ many times.
		
		Plugging the last expression into \eqref{def: eta_n general} with the complete refining property, we obtain
		\begin{align}
			\eta^{\pi, (p)}_n & \le (M+1)^{2p} c^{\frac{3p}{2}} |\pi^n|^{p-1} \Bigg( \sum_{m=0}^{n-1} |\pi^m|^{\frac{1}{p}-\frac{1}{2}} |\pi^{m+1}|^{-\frac{1}{2}} \Big([x]^{(p)}_{\pi^{m+1}}(T)\Big)^{\frac{1}{p}} \Bigg)^p   \nonumber
			\\
			&\le (M+1)^{2p} c^{\frac{3p}{2}} |\pi^n|^{p-1} \Bigg( \sum_{m=0}^{n-1} b^{\frac{1}{2}}|\pi^{m}|^{\frac{1}{p}-1} \Big([x]^{(p)}_{\pi^{m+1}}(T)\Big)^{\frac{1}{p}} \Bigg)^p \nonumber
			\\
			& = (M+1)^{2p} c^{\frac{3p}{2}} b^{\frac{p}{2}} \Bigg( \sum_{m=0}^{n-1} \left(\frac{|\pi^n|}{|\pi^{m}|}\right)^{1-\frac{1}{p}} \Big([x]^{(p)}_{\pi^{m+1}}(T)\Big)^{\frac{1}{p}} \Bigg)^p  \nonumber
			\\
			& \le (M+1)^{2p} c^{\frac{3p}{2}} b^{\frac{p}{2}} \Bigg( \sum_{m=0}^{n-1} \left(1+a\right)^{(m-n)(1-\frac{1}{p})} \Big([x]^{(p)}_{\pi^{m+1}}(T)\Big)^{\frac{1}{p}} \Bigg)^p.      \label{eq : eta general bound}
		\end{align}
		
		We now define an infinite-dimensional matrix $A = (a_{n, m})_{n \ge 0, m \ge 0}$ with entries
		\begin{equation*}
			a_{n, m} := 
			\begin{dcases}
				\Big( 1-(1+a)^{\frac{1}{p}-1} \Big) \times (1+a)^{(m-n)(1-\frac{1}{p})}, \quad &\text{for } m \le n, 
				\\
				\qquad \qquad 0, \qquad &\text{for } m > n,
			\end{dcases}
		\end{equation*}
		and we shall show that the matrix $A$ satisfies properties (i) - (iii) of Lemma~\ref{lem : Silverman-Toeplitz weak}. First, condition (i) is obvious. In order to show (ii), we use the geometric series to derive
		\begin{align*}
			\lim_{n \rightarrow \infty} \sum_{m=0}^{\infty} a_{n, m}
			&= \lim_{n \rightarrow \infty}  \Big( 1-(1+a)^{\frac{1}{p}-1} \Big) \bigg( \sum_{m=0}^{n}  (1+a)^{(m-n)(1-\frac{1}{p})}\bigg)
			\\
			&= \lim _{n \rightarrow \infty} \Big( 1-(1+a)^{\frac{1}{p}-1} \Big) \bigg( \frac{1-(1+a)^{(\frac{1}{p}-1)(n+1)}}{1-(1+a)^{\frac{1}{p}-1}} \bigg)
			\\
			&= \lim _{n \rightarrow \infty} 1-(1+a)^{(\frac{1}{p}-1)(n+1)} = 1.
		\end{align*}
		Condition (iii) is also obvious from (ii); $\sup_{n \ge 0} \sum_{m=0}^{\infty} |a_{n, m}| = 1 < \infty$.
		
		Therefore, we apply Lemma~\ref{lem : Silverman-Toeplitz weak} to the inequality \eqref{eq : eta general bound} to obtain
		\begin{align}   \label{ineq : limsup eta upper bound}
			\limsup_{n \rightarrow \infty} \, \eta^{\pi, (p)}_n
			&\le \frac{(M+1)^{2p} c^{\frac{3p}{2}} b^{\frac{p}{2}}}{\big( 1-(1+a)^{\frac{1}{p}-1} \big)^p} \limsup_{n \rightarrow \infty} \Bigg( \sum_{m=0}^{\infty} a_{n, m} \Big([x]^{(p)}_{\pi^{m+1}}(T)\Big)^{\frac{1}{p}} \Bigg)^p       \nonumber
			\\
			& \le \frac{(M+1)^{2p} c^{\frac{3p}{2}} b^{\frac{p}{2}}}{\big( 1-(1+a)^{\frac{1}{p}-1} \big)^p} \bigg( \limsup_{n \rightarrow \infty} \Big([x]^{(p)}_{\pi^{n}}(T)\Big)^{\frac{1}{p}} \bigg)^p             \nonumber
			\\
			& = \frac{(M+1)^{2p} c^{\frac{3p}{2}} b^{\frac{p}{2}}}{\big( 1-(1+a)^{\frac{1}{p}-1} \big)^p} \limsup_{n \rightarrow \infty} \, [x]^{(p)}_{\pi^{n}}(T).
		\end{align}
		Combining \eqref{ineq : limsup eta upper bound} with the inequality after taking $\limsup$ to \eqref{def: eta general}, yields the result \eqref{ineq : p-th var, eta limsup}.
	\end{proof}
	
\begin{proof} [Proof of Theorem \ref{thm: p-th var general2}]
For fixed $p, x$, and $\pi$ satisfying the conditions of Theorem \ref{thm: p-th var general2}, let us define
\begin{align*}
a_n := \sum_{m=0}^{n-1} \vert \pi^{m} \vert^{\frac{1}{p}-\frac{1}{2}} \bigg( \sum_{k \in I_m} |\theta^{x, \pi}_{m, k}|^p \bigg)^\frac{1}{p}, \qquad \qquad
			b_n := |\pi^n|^{\frac{1}{p}-1}, \qquad \qquad \forall \, n \in \N
		\end{align*}
		such that
		\begin{align*}
			a_{n+1}-a_{n} = \vert \pi^{n} \vert^{\frac{1}{p}-\frac{1}{2}} \bigg( \sum_{k \in I_n} |\theta^{x, \pi}_{n, k}|^p \bigg)^\frac{1}{p},
			\qquad 
			b_{n+1}-b_n = |\pi^{n+1}|^{\frac{1}{p}-1} -|\pi^{n}|^{\frac{1}{p}-1}.
		\end{align*}
		Moreover, from the notation \eqref{def: eta_n general}, we have
		\begin{equation} \label{eq : a_n b_n fractions general}
			\frac{a_n}{b_n} = \big(\eta_n^{\pi, (p)}\big)^{\frac{1}{p}}, \qquad \qquad
			\frac{a_{n+1} - a_n}{b_{n+1} - b_n} 
			= \frac{\vert \pi^{n} \vert^{\frac{1}{p}-\frac{1}{2}}\big( \sum_{k \in I_n} |\theta^{x, \pi}_{n, k}|^p \big)^\frac{1}{p}}{ |\pi^{n+1}|^{\frac{1}{p}-1} -|\pi^{n}|^{\frac{1}{p}-1}} = \frac{\big(\xi^{\pi, (p)}_n\big)^{\frac{1}{p}}}{ \Big(\frac{|\pi^{n+1}|}{|\pi^n|}\Big)^{\frac{1}{p}-1} -1},
		\end{equation}
		and the complete refining property provides the bounds
		\begin{equation}	\label{ineq: a_n b_n diff ratio bounds}
			\frac{\big(\xi^{\pi, (p)}_n\big)^{\frac{1}{p}}}{b^{1-\frac{1}{p}}-1}
			\le 
			\frac{a_{n+1} - a_n}{b_{n+1} - b_n} 
			\le \frac{\big(\xi^{\pi, (p)}_n\big)^{\frac{1}{p}}}{(1+a)^{1-\frac{1}{p}}-1}.
		\end{equation}
		We further define 
		\begin{equation}    \label{def : beta_n}
			\beta_n := \frac{b_{n+1}}{b_n} = \bigg(\frac{|\pi^{n+1}|}{|\pi^n|}\bigg)^{\frac{1}{p}-1} > 1, \qquad \forall \, n \in \mathbb{N},
		\end{equation} 
		then, the limit $\beta := \lim_{n \to \infty} \beta_n = r^{\frac{1}{p}-1} > 1$ exists, thanks to the convergent refining property of $\pi$. Applying \eqref{ineq : limsup} of Lemma \ref{lem: a_n/b_n} with the bounds \eqref{ineq: a_n b_n diff ratio bounds}, \eqref{ineq : limsup eta upper bound} yields
		\begin{align*}
			\limsup_{n \rightarrow \infty} \, \frac{\big( \xi^{\pi, (p)}_n \big)^{\frac{1}{p}}}{b^{1-\frac{1}{p}} -1}
			&\le \frac{\beta}{\beta-1} \limsup_{n \to \infty} \, \big(\eta^{\pi, (p)}_n \big)^{\frac{1}{p}} - \frac{1}{\beta-1} \liminf_{n \to \infty} \, \big(\eta^{\pi, (p)}_n \big)^{\frac{1}{p}}
			\le \frac{\beta}{\beta-1} \limsup_{n \to \infty} \, (\eta^{\pi, (p)}_n)^{\frac{1}{p}} 
			\\
			&\le \bigg(\frac{\beta}{\beta-1}\bigg) \bigg( \frac{(M+1)^{2} c^{\frac{3}{2}} b^{\frac{1}{2}}}{1-(1+a)^{\frac{1}{p}-1}} \bigg) \limsup_{n \to \infty} \, \Big([x]^{(p)}_{\pi^n}6(T)\Big)^{\frac{1}{p}}.
		\end{align*}
		This implies $\limsup_{n \to \infty} \, [x]_{\pi^n}^{(p)}(T) < \infty \Longrightarrow \limsup_{n \to \infty} \, \xi^{\pi, (p)}_n < \infty$.
		
		For the opposite direction, we take $\limsup$ to \eqref{def: eta general}, and use Lemma \ref{lem : liminf and limsup} with \eqref{ineq: a_n b_n diff ratio bounds} to obtain
		\begin{align*}
			\frac{1}{c\big(Mc\sqrt{bM}\big)^p} \limsup_{n \to \infty} \, [x]^{(p)}_{\pi^n}(T) &\le 
			\limsup_{n \rightarrow \infty} \, \eta^{\pi, (p)}_n 
			= \limsup_{n \rightarrow \infty} \, \bigg(\frac{a_n}{b_n}\bigg)^p
			\\
			&\le \limsup_{n \rightarrow \infty} \, \bigg(\frac{a_{n+1}-a_n}{b_{n+1}-b_n}\bigg)^p = \frac{1}{\big((1+a)^{1-\frac{1}{p}}-1\big)^p}\limsup_{n \rightarrow \infty} \, \xi^{\pi, (p)}_n.
		\end{align*}
		This proves the result \eqref{ineq : p-th var, xi limsup}. 
	\end{proof}
	
	\bigskip
	
\section{Isomorphism on \texorpdfstring{$\mathcal{X}^{p}_\pi$}{TEXT}} \label{sec: isomorphism}

In this section, we shall use several function norms and matrix norms, thus, Table 1 lists all the norms with
their definitions for the convenience of readers. In Table \ref{table.norm}, $x$ represents a real-valued continuous function defined on $[0, T]$, and $A$ represents an infinite-dimensional matrix whose $(m, k)$-entry is denoted by $A_{m,k}$.

\begin{table}[h]
\centering
\begin{tabular}{|c|c|}\hline
Function norm & Definition \\ \hline
$\|x\|^{(p)}_\pi$  &  $|x(0)| + \sup_{n \in \N} \, \Big([x]_{\pi^n}^{(p)} (T)\Big)^{\frac{1}{p}}$ \text{ in Definition \eqref{Def : p norm and Xp space}} \\ \hline
$\| x\|_{\infty}$ & $\sup_{t\in [0,T]} |x(t)|$ \\ \hline
$|x|_{C^{0,\alpha}} $ &  $\sup_{s, t \in [0, T], ~ s \neq t} \frac{|x(s) - x(t)|}{|s-t|^{\alpha}}$\\ \hline
$\| x\|_{C^{0,\alpha}} $ &  $\| x\|_\infty + | x|_{C^{0,\alpha}}$ \text{ in \eqref{def : holder norm}}\\ \hhline{|=|=|}
Matrix norm & Definition \\ \hline
$\|A\|_{sup} $ & $\sup_{m,k\geq 0} |A_{m,k}|$\\ \hline
$\|A\|^\alpha_{sup} $ & $\|D^\pi_\alpha A\|_{sup} $ where $D^\pi_\alpha$ is the matrix defined in \eqref{def: D matrix}\\ \hline
$\Vert A \Vert_{p, \infty} $ &  $\sup_{k \ge 0} \big( \sum_{m \ge 0} |A_{m, k}|^p \big)^{\frac{1}{p}}$ in \eqref{def: p infty matrix norm}\\ \hline
$\|A\|_{(p)} $ & $\|(E^\pi A)^{\top}\|_{p, \infty} $ where $E^\pi$ is the matrix defined in \eqref{def: E matrix}\\ \hline
\end{tabular}
\caption{List of norms used in this section}
\label{table.norm}
\end{table}

Recall the space $C^{0, \alpha}([0, T])$ of $\alpha$-H\"older continuous functions with the norm 
\begin{equation}    \label{def : holder norm}
\Vert x \Vert_{C^{0, \alpha}} := \Vert x \Vert_{\infty} + \vert x \vert_{C^{0, \alpha}} \quad \text{with} \quad \Vert x \Vert_{\infty} = \sup_{t \in [0, T]} \vert x(t) \vert \quad \text{and} \quad |x|_{C^{0, \alpha}} := \sup_{\substack{s, t \in [0, T] \\ s \neq t}} \frac{|x(s) - x(t)|}{|s-t|^{\alpha}}.
\end{equation} 
Ciesielski \cite{Ciesielski:isomorphism} proved that the following mapping $T^{\T}_{\alpha}$ is an isomorphism between $C^{0, \alpha}([0, T])$ and the space $\ell^{\infty}(\R)$ of all bounded real sequences, equipped with the supremum norm $\Vert \cdot \Vert_{\infty}$:
	\begin{align*}
		T^{\T}_{\alpha} : C^{0, \alpha}([0, T]) &\xrightarrow{\hspace{10mm}} ~~~ \ell^{\infty}(\mathbb{R})
		\\
		x ~~~ &\xmapsto{\hspace{5mm}} \big\{ 2^{(m+1)(\alpha-\frac{1}{2})}|\theta^{x, \T}_{m, k}|\big\}_{m, k}.
	\end{align*}
	Here, $\theta^{x, \T}_{m, k}$'s are the Schauder coefficients of $x$ along the dyadic partition sequence $\T$, and the double-indexed set $\{ 2^{(m+1)(\alpha-\frac{1}{2})}|\theta^{x, \T}_{m, k}|\}_{m, k}$ can be identified as a real sequence by flattening it. A recent work \citep{fake_fBM} extends this isomorphism to any balanced, complete refining partition sequence $\pi$:
	\begin{align}
		T^{\pi}_{\alpha} : C^{0, \alpha}([0, T]) &\xrightarrow{\hspace{10mm}} ~~~ \ell^{\infty}(\mathbb{R})        \nonumber
		\\
		x ~~~ &\xmapsto{\hspace{5mm}} \left\{ \vert \pi^{m+1} \vert^{\frac{1}{2}-\alpha}|\theta^{x, \pi}_{m, k}|\right\}_{m, k}.   \label{mapping: generalized Ciesielski}
	\end{align}
We may arrange each element of the sequence $\big\{ \vert \pi^{m+1} \vert^{\frac{1}{2}-\alpha}|\theta^{x, \pi}_{m, k}|\big\}_{m, k}$ in a matrix without flattening it. Let us denote $\mathcal{M}$ the space of infinite-dimensional matrices and fix a partition sequence $\pi = (\pi^n)_{n \ge 0}$ of $[0, T]$. For each $m \ge 0$, recall the index set $I_m$ of \eqref{def : index set} corresponding to $\pi$, and consider the subspace
\begin{equation}
\mathcal{M}_{\pi} := \{ A \in \mathcal{M} : A_{m, k} = 0 \quad \text{if } k > |I_m| \} \subset \mathcal M,
	\end{equation}
	composed of infinite-dimensional matrices whose $m$-th row vector can take nonzero values only for the first $|I_m|$ components. We now construct a `Schauder coefficient matrix' $\Theta^{x, \pi}$ in $\mathcal{M}_{\pi}$ to arrange the Schauder coefficients:
\begin{equation*}
(\Theta^{x, \pi})_{m, k} = 
\begin{cases}
	\theta^{x, \pi}_{m, k}, \qquad &\text{if } k \in I_m,
\\
	~ 0, \qquad &\text{otherwise},
\end{cases}
\qquad m \ge 0, \quad k \ge 0.
\end{equation*}
We also define a diagonal matrix $D^{\pi}_{\alpha} \in \mathcal{M}$ with each $(m, m)$-th entry equal to $|\pi^{m+1}|^{\frac{1}{2}-\alpha}$:
\begin{equation}    \label{def: D matrix}
(D^{\pi}_{\alpha})_{m, k} = 
\begin{cases}
|\pi^{m+1}|^{\frac{1}{2}-\alpha}, \qquad &\text{if } m = k,
\\
~~~ 0, \qquad &\text{otherwise},
\end{cases} \qquad m\ge 0,\; k\ge 0.
\end{equation}
From this construction, we have the identity
\begin{equation}    \label{eq : D Theta}
\sup_{m, k} \Big( \vert \pi^{m+1} \vert^{\frac{1}{2}-\alpha}|\theta^{x, \pi}_{m, k}|\Big) = \Vert D^{\pi}_{\alpha} \Theta^{x, \pi} \Vert_{sup},
\end{equation}
where $\Vert A \Vert_{sup} := \sup_{m, k \ge 0} |A_{m, k}|$ is the supremum norm for matrices; in the mapping $T^{\pi}_{\alpha}$ of \eqref{mapping: generalized Ciesielski}, the condition $\big\{ \vert \pi^{m+1} \vert^{\frac{1}{2}-\alpha}|\theta^{x, \pi}_{m, k}|\big\}_{m, k} \in \ell^{\infty}(\R)$ is then equivalent to $\Vert D^{\pi}_{\alpha} \Theta^{x, \pi} \Vert_{sup} < \infty$.
	
We now restate the isomorphism in \eqref{mapping: generalized Ciesielski} along any balanced and complete refining partition sequence.
	
	\begin{proposition} \label{prop : Isomorphism}
		For any balanced, complete refining partition sequence $\pi$ and any $\alpha \in (0, 1)$, the mapping 
		\begin{align}
			T^{\pi}_{\alpha} : \Big( C^{0, \alpha}([0, T]), \, \Vert \cdot \Vert_{C^{0, \alpha}} \Big) &\xrightarrow{\hspace{10mm}} \Big( \mathcal{M}^{\alpha}_{\pi}, \, \Vert \cdot \Vert^{\alpha}_{sup} \Big)        \nonumber
			\\
			x ~~~~~~~~~~ &\xmapsto{\hspace{10mm}}  ~~~~~ \Theta^{x, \pi}  \label{mapping: generalized Schauder matrix}
		\end{align}
		is an isomorphism, where
		\begin{align*}
			\mathcal{M}^{\alpha}_{\pi} := \{ A \in \mathcal{M}_{\pi} : \Vert A \Vert^{\alpha}_{sup} < \infty\}, \qquad \Vert A \Vert^{\alpha}_{sup} := \Vert D^{\pi}_{\alpha} A \Vert_{sup}.
		\end{align*}
		Moreover, we have the following bounds for the operator norms:
		\begin{equation}
			\Vert T^{\pi}_{\alpha} \Vert_{op} \le 2(\sqrt{c})^3, \qquad 
			\Vert (T^{\pi}_{\alpha})^{-1} \Vert_{op} \le \max \Big( 2M\sqrt{c}K^{\alpha}_1 + 2MK^{\alpha}_2, \, MK^{\alpha}_2 |\pi^1|^{\alpha} \Big),
		\end{equation}
		where $K^{\alpha}_1 := \frac{1}{1-(1+a)^{\alpha-1}}$ and $K^{\alpha}_2 := \frac{1}{1-(1+a)^{-\alpha}}$ with the constants $a, c, M$ in Remark \ref{rem: constants}.
	\end{proposition}
	
	\begin{proof}[Proof of Proposition \ref{prop : Isomorphism}]
		From \cite[Theorem 3.4]{fake_fBM} and the identity \eqref{eq : D Theta}, it is easy to show that the mapping $T^{\pi}_{\alpha}$ is bijective. We note that the notation $\Vert \cdot \Vert_{C^{\alpha}([0, T])}$ in the bounds \cite[Equation (3.2)]{fake_fBM} represents the H\"older semi-norm ($\vert \cdot \vert_{C^{0, \alpha}}$ in \eqref{def : holder norm} of this paper).
		
The bound for operator norm $\Vert T^{\pi}_{\alpha} \Vert_{op}$ is also straightforward from \cite[Theorem 3.4]{fake_fBM} and \eqref{eq : D Theta}:
\begin{equation*}
	\Vert \Theta^{x, \pi} \Vert^{\alpha}_{sup} = \sup_{m, k} \Big( \vert \pi^{m+1} \vert^{\frac{1}{2}-\alpha}|\theta^{x, \pi}_{m, k}|\Big) \le 2(\sqrt{c})^3 \vert x \vert_{C^{0, \alpha}} \le 2(\sqrt{c})^3 \Vert x \Vert_{C^{0, \alpha}}.
\end{equation*}
		The same theorem also yields the inequality
		\begin{equation}    \label{ineq : semi norm bound}
			\vert x \vert_{C^{0, \alpha}} \le (2M\sqrt{c}K^{\alpha}_1 + 2MK^{\alpha}_2) \Vert \Theta^{x, \pi} \Vert^{\alpha}_{sup}.
		\end{equation}
		Furthermore, we can derive that
		\begin{align*}
			\Vert x \Vert_{\infty} 
			&\le \sup_{t \in [0, T]} \bigg( \sum_{m = 0}^{\infty} \sum_{k \in I_m} \vert \theta^{x, \pi}_{m, k} \vert \vert e^{\pi}_{m, k}(t) \vert \bigg)
			\le M \sum_{m = 0}^{\infty} \Big( \sup_{k \in I_m} |\theta^{x, \pi}_{m, k}| \Big) |\pi^{m+1}|^{\frac{1}{2}}
			\\
			&\le M \Big( \sum_{m = 0}^{\infty} |\pi^{m+1}|^{\alpha} \Big) \bigg( \sup_{m, k} \Big(|\theta^{x, \pi}_{m, k}||\pi^{m+1}|^{\frac{1}{2}-\alpha} \Big) \bigg)
			\le MK^{\alpha}_2 |\pi^1|^{\alpha} \Vert \Theta^{x, \pi} \Vert^{\alpha}_{sup}.
		\end{align*}
		Here, the second inequality and the last inequality follow from  \cite[bound (2.4) and Lemma 3.2]{fake_fBM}, respectively. Combining this with \eqref{ineq : semi norm bound} yields the bound for $\Vert (T^{\pi}_{\alpha})^{-1} \Vert_{op}$.
	\end{proof}
	
	Let us fix $x \in C^{0, \alpha}([0, T])$ and $\pi \in \Pi([0, T])$, and recall from Theorem \ref{thm : variation index bounds} that $x$ belongs to $\mathcal{X}^{q}_{\pi}$ for some $q \in [1, \frac{1}{\alpha}]$. In what follows, we shall characterize such functions $x \in C^{0, \alpha}([0, T]) \cap \mathcal{X}^{q}_{\pi}$ in terms of its Schauder coefficients.
	
	We now fix $p > 1$ and define a diagonal matrix $E^{\pi}$ in $\mathcal{M}$ such that every $(m, m)$-th entry is equal to $|\pi^m|^{\frac{1}{2}}$:
    \begin{equation}    \label{def: E matrix}
        (E^{\pi})_{m, k} := 
        \begin{cases}
		  |\pi^{m}|^{\frac{1}{2}}, \qquad &\text{if } m = k,
		  \\
		  ~~~ 0, \qquad &\text{otherwise}, 
		\end{cases}\qquad m\ge 0,\; k\ge 0.
    \end{equation}
    With the matrix norm
	\begin{equation}   \label{def: p infty matrix norm}
	\Vert A \Vert_{p, \infty} := \sup_{k \ge 0} \Big( \sum_{m \ge 0} |A_{m, k}|^p \Big)^{\frac{1}{p}}, \qquad \text{for any } p > 1,
	\end{equation}
	we define
	\begin{equation}    \label{def : (p) matrix norm}
		\mathcal{M}^{(p)}_{\pi} := \{ A \in \mathcal{M}_{\pi} : \Vert A \Vert_{(p)} < \infty\}, \qquad \text{where} \qquad \Vert A \Vert_{(p)} := \Vert (E^{\pi} A)^{\top} \Vert_{p, \infty}.
	\end{equation}
	Recalling the definition \eqref{def : xi general}, we obtain the identity from \eqref{def : (p) matrix norm}
	\begin{equation}    \label{eq : (p) matrix norm and sup xi}
		\Vert \Theta^{x, \pi} \Vert_{(p)} = \Vert (E^{\pi}\Theta^{x, \pi})^{\top} \Vert_{p, \infty} = \sup_{n \ge 0} \big( \xi^{\pi, (p)}_n \big)^{\frac{1}{p}}.
	\end{equation}
	Therefore, the condition \eqref{ineq : p-th var, xi limsup} of Theorem \ref{thm: p-th var general2} is also equivalent to $\Vert \Theta^{x, \pi} \Vert_{(p)} < \infty$. We are now ready to provide the following results regarding the intersection space $C^{0, \alpha}([0, T]) \cap \mathcal{X}^{p}_{\pi}$.
	
	\begin{proposition}\label{prop. intersection banach space}
		For any $\alpha \in (0, 1)$, $p \in (1, \frac{1}{\alpha}]$, and $\pi \in \Pi([0, T])$, the space $\big( C^{0, \alpha}([0, T]) \cap \mathcal{X}^{p}_{\pi}, \, \Vert \cdot \Vert_{C^{0, \alpha}} + \Vert \cdot \Vert^{(p)}_{\pi} \big)$ is a Banach space.
	\end{proposition}
	
	\begin{proof}[Proof of Proposition \ref{prop. intersection banach space}]
		Since $\big(C^{0, \alpha}([0, T]), \, \Vert \cdot \Vert_{C^{0, \alpha}} \big)$ and $\big( \mathcal{X}^{p}_{\pi}, \, \Vert \cdot \Vert^{(p)}_{\pi} \big)$ are Banach spaces (Proposition \ref{prop : Banach space}), it is obvious that $\Vert \cdot \Vert_{C^{0, \alpha}} + \Vert \cdot \Vert^{(p)}_{\pi}$ is a norm in the intersection space, and it is enough to show the completeness of $C^{0, \alpha}([0, T]) \cap \mathcal{X}^{p}_{\pi}$. Fix any Cauchy sequence $(x_{\ell})_{\ell \in \mathbb{N}} \in C^{0, \alpha}([0, T]) \cap \mathcal{X}^{p}_{\pi}$ in $\Vert \cdot \Vert_{C^{0, \alpha}} + \Vert \cdot \Vert^{(p)}_{\pi}$-norm. Then, $(x_{\ell})_{\ell \in \mathbb{N}}$ is also Cauchy in $\Vert \cdot \Vert_{C^{0, \alpha}}$-norm, thus it has a limit $x \in C^{0, \alpha}([0, T])$ such that $\Vert x_{\ell} - x \Vert_{C^{0, \alpha}} \to 0$ as $\ell \to \infty$; in particular, $\{x_{\ell}(t)\}_{\ell \in\mathbb{N}}$ is a Cauchy sequence in $\mathbb{R}$, and $x_{\ell}(t) \to x(t)$ as $\ell \to \infty$ for each $t \in [0, T]$. Moreover, since $\{x_{\ell}\}_{\ell \in\mathbb{N}}$ is also a Cauchy sequence in $\Vert \cdot \Vert^{(p)}_{\pi}$-norm, there exists a limit $\tilde{x} \in \mathcal{X}^{p}_{\pi}$ such that $\Vert x_{\ell} - \tilde{x} \Vert^{(p)}_{\pi} \to 0$ as $\ell \to \infty$. As in the proof of Proposition \ref{prop : Banach space}, we have $\lim_{\ell \to \infty} x_{\ell}(t^n_j) = \tilde{x}(t^n_j) = x(t^n_j)$ for every partition point $t^n_j$ of $P := \bigcup_{n \ge 0} \pi^n$. In other words, $x$ and $\tilde{x}$ coincide on the dense set $P$, thus the unique continuous extension of $\tilde{x}$ must be $x$, thus $(x_{\ell})_{\ell\in\mathbb{N}}$ converges to $x \in C^{0, \alpha}([0, T]) \cap \mathcal{X}^{p}_{\pi}$ in $\Vert \cdot \Vert_{C^{0, \alpha}} + \Vert \cdot \Vert^{(p)}_{\pi}$-norm.
	\end{proof}
	
	In addition to Ciesielski's isomorphism, we have the following isomorphism from the intersection space.
	\begin{theorem} [Isomorphism on the Banach space $\mathcal{X}^{p}_{\pi}$] \label{thm : isomorphism}
		For any $\alpha \in (0, 1)$, $p \in (1, \frac{1}{\alpha}]$, and a balanced, convergent refining partition sequence $\pi$, the mapping
		\begin{align}
			T^{\pi}_{\alpha, (p)} : \Big( C^{0, \alpha}([0, T]) \cap \mathcal{X}^{p}_{\pi}, \, \Vert \cdot \Vert_{C^{0, \alpha}} + \Vert \cdot \Vert^{(p)}_{\pi} \Big) &\xrightarrow{\hspace{10mm}} \Big( \mathcal{M}^{\alpha}_{\pi} \cap \mathcal{M}^{(p)}_{\pi}, \, \Vert \cdot \Vert^{\alpha}_{sup} + \Vert \cdot \Vert_{(p)} \Big)        \nonumber
			\\
			x ~~~~~~~~~~ &\xmapsto{\hspace{10mm}}  ~~~~~ \Theta^{x, \pi}  \label{mapping: isomorphism}
		\end{align}
		is an isomorphism. Furthermore, we have the following bounds for the operator norms:
		\begin{align}
			\Vert T^{\pi}_{\alpha, (p)} \Vert_{op} &\le \max \bigg( 2(\sqrt{c})^3, \, \frac{(M+1)^{2}c^{\frac{3}{2}}b^{\frac{3}{2}-p}}{\Big((1+a)^{1-\frac{1}{p}}-1\Big)^{\frac{1}{p}}} \bigg),  \label{ineq : T alpha p bound}
			\\
			\Vert (T^{\pi}_{\alpha, (p)})^{-1} \Vert_{op} &\le 1 + \max \Big( 2M\sqrt{c}K^{\alpha}_1 + 2MK^{\alpha}_2, \, MK^{\alpha}_2 |\pi^1|^{\alpha} \Big) + \frac{c^{\frac{1}{p}} (Mc\sqrt{bM})}{(1+a)^{1-\frac{1}{p}}-1}.   \label{ineq : T alpha p inverse bound}
		\end{align}
	\end{theorem}
	
	\begin{proof}[Proof of Theorem \ref{thm : isomorphism}] 
		We shall prove the result in the following parts.
		
		\noindent \textbf{Part 1:} For any $x \in C^{0, \alpha}([0, T]) \cap \mathcal{X}^{p}_{\pi}$, we shall prove $T^{\pi}_{\alpha, (p)} (x) \in \mathcal{M}^{\alpha}_{\pi} \cap \mathcal{M}^{(p)}_{\pi}$. \\
		We fix $x \in C^{0, \alpha}([0, T]) \cap \mathcal{X}^{p}_{\pi}$. Proposition \ref{prop : Isomorphism} proves $\Theta^{x, \pi} \in \mathcal{M}^{\alpha}_{\pi}$, thus we need to show $\Theta^{x, \pi} \in \mathcal{M}^{(p)}_{\pi}$, which is equivalent to $\sup_{n \ge 0} \big( \xi^{(p)}_n \big) < \infty$ from \eqref{eq : (p) matrix norm and sup xi}.
		
		Recalling the inequality \eqref{eq : eta general bound} and computing the geometric series, we have for each $n \ge 0$
\begin{align*}			 \eta^{\pi, (p)}_{\pi^n}
&	\le (M+1)^{2p}c^{\frac{3p}{2}}b^{\frac{p}{2}} \Big( \Vert x \Vert^{(p)}_{\pi} \Big)^p \bigg( \sum_{m=0}^{n-1} \left(1+a\right)^{(m-n)(1-\frac{1}{p})} \bigg)^p
			\\
			& = (M+1)^{2p}c^{\frac{3p}{2}}b^{\frac{p}{2}} \Big( \Vert x \Vert^{(p)}_{\pi} \Big)^p \bigg( \frac{1-(1+a)^{-n(1-\frac{1}{p})}}{(1+a)^{1-\frac{1}{p}}-1}\bigg)
			\le \frac{(M+1)^{2p}c^{\frac{3p}{2}}b^{\frac{p}{2}}}{(1+a)^{1-\frac{1}{p}}-1} \Big( \Vert x \Vert^{(p)}_{\pi} \Big)^p.
		\end{align*}
		Furthermore, recalling the notations \eqref{eq : a_n b_n fractions general} and \eqref{def : beta_n} with the identity \eqref{eq : a_n b_n identity}, we derive
		\begin{align*}
			\big( \xi^{\pi, (p)}_n \big)^{\frac{1}{p}} = (\beta_n-1)\frac{a_{n+1} - a_n}{b_{n+1} - b_n}
			= \beta_n \frac{a_{n+1}}{b_{n+1}} - \frac{a_n}{b_n}
			\le \beta_n \frac{a_{n+1}}{b_{n+1}}
			\le b^{1-\frac{1}{p}} \Big( \eta^{\pi, (p)}_{n+1} \Big)^{\frac{1}{p}}.
		\end{align*}
		Here, the last inequality uses the fact that $\beta_n$ has an upper bound $b^{1-\frac{1}{p}}$ from the complete refining property.
		
		Combining the last two inequalities, we obtain for each $n \ge 0$
		\begin{equation}    \label{ineq : xi upper bound}
			\big( \xi^{\pi, (p)}_n \big)^{\frac{1}{p}} \le \frac{(M+1)^{2}c^{\frac{3}{2}}b^{\frac{3}{2}-p}}{\Big((1+a)^{1-\frac{1}{p}}-1\Big)^{\frac{1}{p}}} \Vert x \Vert^{(p)}_{\pi}.
		\end{equation}
		Since $x \in \mathcal{X}^{p}_{\pi}$, we have $\sup_{n \ge 0} \big( \xi^{(p)}_n \big) < \infty$, which shows $\Theta^{x, \pi} \in \mathcal{M}^{(p)}_{\pi}$.
		
		\medskip
		
		\noindent \textbf{Part 2:} For any $\Theta \in \mathcal{M}^{\alpha}_{\pi} \cap \mathcal{M}^{(p)}_{\pi}$, we shall prove $(T^{\pi}_{\alpha, (p)})^{-1}\Theta \in C^{0, \alpha}([0, T]) \cap \mathcal{X}^{p}_{\pi}$. \\
		We fix $\Theta \in \mathcal{M}^{\alpha}_{\pi} \cap \mathcal{M}^{(p)}_{\pi}$. Using the entries $\Theta_{m, k}$ of $\Theta$ as Schauder coefficients along $\pi$, we can construct an $\alpha$-H\"older continuous function $x$ from Proposition \ref{prop : Isomorphism}. The identity \eqref{eq : (p) matrix norm and sup xi} with Corollary \ref{cor : new def variation index general2} and \eqref{def. alter variation index2} imply $x \in \mathcal{X}^{p}_{\pi}$.
		
		\medskip
		
		\noindent \textbf{Part 3:} We shall prove that the mapping $T^{\pi}_{\alpha, (p)}$ is bounded.
		
		For any $x \in C^{0, \alpha}([0, T]) \cap \mathcal{X}^{p}_{\pi}$, consider $\Theta^{x, \pi} = T^{\pi}_{\alpha, (p)}x$. From \eqref{eq : (p) matrix norm and sup xi} and \eqref{ineq : xi upper bound}, we have
		\begin{equation*}
			\Vert \Theta^{x, \pi} \Vert_{(p)} \le \frac{(M+1)^{2}c^{\frac{3}{2}}b^{\frac{3}{2}-p}}{\Big((1+a)^{1-\frac{1}{p}}-1\Big)^{\frac{1}{p}}} \Vert x \Vert^{(p)}_{\pi}.
		\end{equation*}
		Moreover, from Proposition \ref{prop : Isomorphism}, we have $\Vert \Theta^{x, \pi} \Vert^{\alpha}_{sup} \le 2(\sqrt{c})^3 \Vert x \Vert_{C^{0, \alpha}}$.
		Combining the two bounds concludes \eqref{ineq : T alpha p bound}.
		
		\medskip
		
		\noindent \textbf{Part 4:} We shall prove that the inverse mapping $(T^{\pi}_{\alpha, (p)})^{-1}$ is bounded.
		
		For any $\Theta \in \mathcal{M}^{\alpha}_{\pi} \cap \mathcal{M}^{(p)}_{\pi}$, we write $x = (T^{\pi}_{\alpha, (p)})^{-1}\Theta$ and consider its Schauder coefficients $\{\theta^{x, \pi}_{m, k} = \Theta_{m, k}\}_{m, k}$.
		Recalling the inequality \eqref{def: eta general} and the notation \eqref{def : xi general}, we obtain for any $n \ge 0$
		\begin{align*}
			[x]_{\pi^n}^{(p)}(T) 
			& \le \Big(Mc\sqrt{bM}|\pi^n|\Big)^p \Bigg( \sum_{m=0}^{n-1} \vert \pi^{m} \vert^{-\frac{1}{2}} \bigg( \frac{c|\pi^{m}|}{|\pi^n|}\bigg)^\frac{1}{p} \bigg( \sum_{k, i} |\theta^{x, \pi}_{m, k, i}|^p \bigg)^\frac{1}{p} \Bigg)^{p}
			\\
			& \le \Big(Mc\sqrt{bM}|\pi^n|\Big)^p \Bigg( \sum_{m=0}^{n-1} \vert \pi^{m} \vert^{-\frac{1}{2}} \bigg( \frac{c|\pi^{m}|}{|\pi^n|}\bigg)^\frac{1}{p} |\pi^m|^{-\frac{1}{2}} (\xi^{\pi, (p)}_m)^{\frac{1}{p}} \Bigg)^{p}
			\\
			& = c \Big(Mc\sqrt{bM}\Big)^p |\pi^n|^{p-1} \Bigg( \sum_{m=0}^{n-1} \vert \pi^{m} \vert^{\frac{1}{p}-1} \bigg)^p \Big(\sup_{m \ge 0} \, \xi^{\pi, (p)}_m\Big).
		\end{align*}
		From the complete refining property and computing the geometric series, we have for each $n \ge 0$
		\begin{align*}
			\sum_{m=0}^{n-1} |\pi^m|^{\frac{1}{p}-1}
			&\le |\pi^n|^{\frac{1}{p}-1} \sum_{m=0}^{n-1} (1+a)^{(\frac{1}{p}-1)(n-m)}
			\\
			&= |\pi^n|^{\frac{1}{p}-1} (1+a)^{\frac{1}{p}-1}\frac{1-(1+a)^{(\frac{1}{p}-1)n}}{1-(1+a)^{\frac{1}{p}-1}}
			\le |\pi^n|^{\frac{1}{p}-1}
			\frac{(1+a)^{\frac{1}{p}-1}}{1-(1+a)^{\frac{1}{p}-1}}
			= \frac{|\pi^n|^{\frac{1}{p}-1}}{(1+a)^{1-\frac{1}{p}}-1}.
		\end{align*}
		Combining the last two inequalities,
		\begin{align*}
			[x]_{\pi^n}^{(p)}(T)
			\le c \Big(Mc\sqrt{bM}\Big)^p |\pi^n|^{p-1} \Bigg( \frac{|\pi^n|^{\frac{1}{p}-1}}{(1+a)^{1-\frac{1}{p}}-1} \Bigg)^p \Big(\sup_{m \ge 0} \, \xi^{\pi, (p)}_m\Big)
			= \frac{c \Big(Mc\sqrt{bM}\Big)^p}{\Big((1+a)^{1-\frac{1}{p}}-1\Big)^p} \Big(\sup_{m \ge 0} \, \xi^{\pi, (p)}_m\Big).
		\end{align*}
		Moreover, thanks to \eqref{eq : (p) matrix norm and sup xi}, we have
		\begin{equation*}
			\Vert x \Vert^{(p)}_{\pi} 
			\le |x(0)| + \frac{c^{\frac{1}{p}} (Mc\sqrt{bM})}{(1+a)^{1-\frac{1}{p}}-1} \Big(\sup_{m \ge 0} \, \xi^{\pi, (p)}_m\Big)^{\frac{1}{p}}
			= |x(0)| + \frac{c^{\frac{1}{p}} (Mc\sqrt{bM})}{(1+a)^{1-\frac{1}{p}}-1} \Vert \Theta^{x, \pi} \Vert_{(p)}.
		\end{equation*}
		Also, Proposition \ref{prop : Isomorphism} yields a bound $\Vert x \Vert_{C^{0, \alpha}} \le \max \Big( 2M\sqrt{c}K^{\alpha}_1 + 2MK^{\alpha}_2, \, MK^{\alpha}_2 |\pi^1|^{\alpha} \Big) \Vert \Theta \Vert^{\alpha}_{sup}$. Combining these bounds proves \eqref{ineq : T alpha p inverse bound}.
	\end{proof}
	
\begin{remark}
		From Proposition \ref{prop : Isomorphism} and Theorem \ref{thm : isomorphism}, one may expect that the following mapping would also be an isomorphism:
		\begin{align*}
			T^{\pi}_{(p)} : \Big( \mathcal{X}^{p}_{\pi}, \, \Vert \cdot \Vert^{(p)}_{\pi} \Big) &\xrightarrow{\hspace{10mm}} \Big( \mathcal{M}^{(p)}_{\pi}, \, \Vert \cdot \Vert_{(p)} \Big)        \nonumber
			\\
			x ~~~~~~~~~~ &\xmapsto{\hspace{10mm}}  ~~~~~ \Theta^{x, \pi}.
		\end{align*}
		However, this is not an isomorphism, since $x \in \mathcal X ^{(p)}_\pi$ is a subclass of continuous functions, and the continuity is not guaranteed without additional conditions if one constructs a function from Schauder coefficients. In the following, we provide an example of function $x$ constructed from a given Schauder matrix $\Theta \in \mathcal{M}^{(2)}_{\pi}$, satisfying the condition $\Vert x \Vert^{(2)}_{\pi} < \infty$, but $x \notin C^0([0,T],\R)$.
		
		Let us consider the dyadic partition sequence $\mathbb{T}$ on a unit interval $[0, 1]$ and a matrix $\Theta \in \mathcal{M}$ such that for each $m \ge 0$ the components of $m$-th row are given by $\Theta_{m, 0} = 2^{\frac{m}{2}}$ and $\Theta_{m, k} = 0$  for all $k \ge 1$. Then, it is easy to verify that $\Vert \Theta \Vert_{(2)} = \Vert (E^{\T} \Theta)^{\top} \Vert_{2, \infty} < \infty$. We now construct a function $x(\cdot) := \sum_{m=0}^{\infty}\sum_{k \in I_m} \Theta_{m, k}e^{\T}_{m, k}(\cdot)$ on $[0, 1]$. It turns out that $x$ is not continuous at $0$; we take $t_n = 2^{-n}$ for each $n \in \N$, then we have
		\begin{equation*}
			x(t_n) = \sum_{m=0}^{n-1} \Theta_{m,0} e^{\T}_{m,0}(t_n) = \sum_{m=0}^{n-1} 2^{\frac{m}{2}} 2^{\frac{m}{2}} t_n = 2^{-n}\sum_{m=0}^{n-1}2^m = 1-2^{-n},
		\end{equation*}
		thus $0 = x(0) = x(\lim_{n \to \infty} t_n) \neq \lim_{n \to \infty} x(t_n) = 1$, so $x\notin C^0([0,1],\R)$.
\end{remark}

\renewcommand{\bibname}{References}
\bibliography{pathwise1}
\bibliographystyle{apalike}

\end{document}